# SCHEDULING CONTROL FOR QUEUEING SYSTEMS WITH MANY SERVERS: ASYMPTOTIC OPTIMALITY IN HEAVY TRAFFIC[1]

By Rami Atar

*Technion–Israel Institute of Technology*

A multiclass queueing system is considered, with heterogeneous service stations, each consisting of many servers with identical capabilities. An optimal control problem is formulated, where the control corresponds to scheduling and routing, and the cost is a cumulative discounted functional of the system's state. We examine two versions of the problem: "nonpreemptive," where service is uninterruptible, and "preemptive," where service to a customer can be interrupted and then resumed, possibly at a different station. We study the problem in the asymptotic heavy traffic regime proposed by Halfin and Whitt, in which the arrival rates and the number of servers at each station grow without bound. The two versions of the problem are not, in general, asymptotically equivalent in this regime, with the preemptive version showing an asymptotic behavior that is, in a sense, much simpler. Under appropriate assumptions on the structure of the system we show: (i) The value function for the preemptive problem converges to $V$, the value of a related diffusion control problem. (ii) The two versions of the problem are asymptotically equivalent, and in particular nonpreemptive policies can be constructed that asymptotically achieve the value $V$. The construction of these policies is based on a Hamilton–Jacobi–Bellman equation associated with $V$.

## 1. Introduction.

We consider a queueing system with a fixed number of customer classes and a fixed number of service stations, each consisting of many servers with identical capabilities. Only some stations can offer service to each class, and the service rates depend on the class and on the station. Customers have a single service requirement and leave the system after it is met. They may also abandon while waiting to be served [see Figure 1(a)].

---

Received December 2003; revised March 2005.

[1]Supported in part by Israel Science Foundation Grant 126/02.

*AMS 2000 subject classifications.* 60K25, 68M20, 90B22, 90B36, 60F05, 49L20.

*Key words and phrases.* Multiclass queueing systems, heavy traffic, scheduling and routing, asymptotically optimal controls, Hamilton–Jacobi–Bellman equation.

Arrivals are modeled as renewal processes, while service and abandonments occur according to exponential clocks. Scheduling of jobs, which amounts to selecting what customer each server serves at each time, is considered as control. The goal is to minimize an expected cumulative discounted functional of general performance criteria, including queue lengths of different classes, number of customers in each of the stations and number of servers that are idle at each station. This paper analyzes the scheduling control (SC) problem in a central limit theorem (CLT) asymptotic regime. In particular, as proposed first by Halfin and Whitt [8], system parameters are rescaled in such a way that the arrival rates and the number of servers at each station are proportional to a parameter $n$, and the system is kept in heavy traffic as $n$ becomes large. See [7] for motivation for this model and parametric regime, in particular in relation to call center models.

One important feature of a CLT approach is that the non-Markovian structure of the SC problem becomes Markovian as weak limits are taken formally. In particular, a diffusion control (DC) problem arises in the limit, that is accessible by control-theoretic tools. In [2] we studied the dynamic programming partial differential equation (PDE) of Hamilton–Jacobi–Bellman (HJB) type for the DC problem and characterized the value function as its unique solution. It is natural to ask whether the value function for the SC problem converges to that of the DC problem.

In fact, this question must be formulated more carefully. One could consider two versions of the problem (that formally yield two different DC problems): one in which only nonpreemptive scheduling is possible, and one where it is allowed to use preemptive policies, where service to a customer can be interrupted and then resumed, possibly at a different station (it is the DC problem associated with preemptive scheduling that was analyzed in [2]). Interestingly, these two versions show identical asymptotic behavior in some cases of the problem, and they show different behavior in other cases: In the single-station case, it was proved in [3] that, in both versions of the problem, the value function for the SC problem converges to that of the same DC problem. These results should be regarded as a simplification that shows in the limit, because the DC problem that is associated formally with nonpreemptive scheduling lies in higher dimension. On the other hand, we will present a heuristic argument suggesting that for a particular two-station system one does not expect the value function in both versions of the problem to have identical limits. Hence one would like to understand when the two versions of the problem behave the same asymptotically.

To indicate what prevents the two-station example from undergoing the same simplification, we need the following terminology. A class-station pair $(i,j)$ is said to be an *activity* if servers of station $j$ are capable of serving class-$i$ customers [in Figure 1, $(1,A)$ is an activity, and $(2,A)$ is not]. A distinction between two kinds of activities is made in terms of an underlying static *fluid model* (introduced in Section 2): The fluid model has a positive fraction of the customer population in some activities and a zero fraction in others. Such activities are referred to as *basic* and, respectively, *nonbasic* activities (cf. [10]). Roughly, the population of customers in basic and nonbasic activities in the *queueing model* is typically $O(n)$ and, respectively, $O(n^{1/2})$. Now, while under a heavy traffic assumption single-station systems have only basic activities; the two-station example alluded to above has a nonbasic activity. As explained in detail in Section 2.7, the simplification will not occur in the presence of an activity having only $O(n^{1/2})$ customers. Aiming at extending the results, that show simplification of the single-station nonpreemptive problem to the multistation setting, we are led to impose the assumption that, in fact, *all activities are basic*. In this paper we find that under this condition, a simplification indeed takes place.

It follows from a result of Williams [14] that if all activities are basic, then the graph, having classes and stations as vertices and activities as edges, is necessarily a tree [as, e.g., in Figure 1(b)]. We refer to this structural condition by saying that the system is *treelike*. As in [2], the results of this paper will therefore apply to treelike systems only.

The main results are the following. Under appropriate assumptions, the preemptive SC problem's value converges to the value of the DC problem. Moreover, preemptive and nonpreemptive scheduling control policies (SCPs) are constructed, that asymptotically achieve the DC problem's value. These results establish the validity of the DC problem as the correct asymptotic description of the problem, as well as the asymptotic equivalence of the preemptive and nonpreemptive problems.

In a preemptive SC problem, one can consider the vector representing the number of customers of each class present in the system as "state," and the vector representing the number of customers at each activity as "control" (at least in the Markovian case, i.e., under the assumption that arrivals are Poisson). In fact, "state" and "control" for the DC problem stand for the formal weak limits of precisely these quantities. The results of [2] regarding the DC problem assert that there exists an optimal Markov control policy, sometimes referred to as an optimal feedback function. Namely, there is a function mapping state to control, that can be used to define the control process as a feedback from the state process, in such a way that the DC problem's value is achieved. In this paper, the preemptive SC problem is solved, roughly, by dynamically setting the population in the activities equal to the image of the system's state under the feedback function for the DC problem. In the nonpreemptive version of the problem one does not have direct control over the population in the activities, and therefore a scheme as above cannot be applied. Instead one must affect the population in the activities indirectly by appropriate decisions of job allocations. To this end, a kind of tracking mechanism is introduced, that roughly works as follows. One continuously computes the desired value of the population in all activities as the image of the state under the feedback function. One declares activities for which the actual population is greater than desired as "overpopulated." Overpopulated activities are then blocked in the sense that new jobs are only assigned through activities that are not overpopulated. As a result, the population decreases rapidly in overpopulated activities and increases rapidly in other activities. The main difficulty dealt with in this paper consists of estimating the difference between the desired and actual values so as to show that tracking is maintained. The proof requires two different types of arguments: one has to do with tightness of the processes involved, and the other regards estimates on their large time behavior. Although the methods of [3] are used in various places in the proof, the main difficulty alluded to above is far more complicated in the multistation case.

We point out that in the more "standard" CLT regime, where the capacity of each server is scaled up rather than the number of servers, a diffusion model for a similar queueing system has been studied by Harrison and López [10]. A reduction, quite different from the one alluded to above, takes place in this regime, that in fact makes the diffusion model one-dimensional. Recently, a number of works have developed policies that are asymptotically optimal in this regime, for models as in [10] as well as more general ones [1, 4, 13].

The organization of the paper is as follows. In Section 2 we introduce the probabilistic queueing model, the scaling and some assumptions regarding work conservation. We then describe the DC problem and the HJB equation, propose SCPs for the queueing model that are based on the HJB equation and state the main result regarding asymptotic optimality of these policies. Finally we discuss the roles of the main assumptions. Section 3 contains the proofs.

*Notation and terminology.* Vectors in $\mathbb{R}^k$ are considered as column vectors. For $x \in \mathbb{R}^k$ (resp., $x \in \mathbb{R}^{k\times l}$) let $\|x\| = \sum_i |x_i|$ (resp., $\sum_{i,j} |x_{i,j}|$). For two column vectors $v, u$, $v \cdot u$ denotes their scalar product. The symbols $e_i$ denote the unit coordinate vectors and $e = (1, \dots, 1)'$. The dimension of $e$ may change from one expression to another, and for example, $e \cdot a + e \cdot b = \sum_i a_i + \sum_j b_j$ even if $a$ and $b$ are of different dimension. Write $\mathbb{N} = \{1, 2, \dots\}$, $\mathbb{Z}_+ = \{0, 1, 2, \dots\}$, $\mathbb{R}_+ = [0, \infty)$. The class of twice continuously differentiable functions on $\mathbb{R}^I$ is denoted by $C^2$. The class of bounded functions in $C^2$ is denoted by $C^2_b$, and $C_{\mathrm{pol}}$ denotes the class of functions $f$ in $C^2$ satisfying a polynomial growth condition: there are constants $c$ and $r$ such that $|f(x)| \le c(1 + \|x\|^r)$, $x \in \mathbb{R}^I$. For $E$ a metric space, we denote by $\mathbb{D}(E)$ the space of all cadlag functions (i.e., right-continuous and having

left limits) from $\mathbb{R}_+$ to $E$. We endow $\mathbb{D}(E)$ with the usual Skorohod topology (cf. [5]). If $X^n$, $n \in \mathbb{N}$, and $X$ are processes with sample paths in $\mathbb{D}(E)$, we write $X^n \Rightarrow X$ to denote weak convergence of the measures induced by $X^n$ [on $\mathbb{D}(E)$] to the measure induced by $X$. For a collection $\mathcal{A}$ of random variables, $\sigma\{\mathcal{A}\}$ denotes the sigma-field generated by this collection. If $X$ is an $\mathbb{R}^k$- or an $\mathbb{R}^{k\times l}$-valued process (or function on $\mathbb{R}_+$) and $0 \le s \le t < \infty$, $\|X\|^*_{s,t} = \sup_{s\le u\le t} \|X(u)\|$, and if $X$ takes real values, $|X|^*_{s,t} = \sup_{s\le u\le t} |X(u)|$. Also, $\|X\|^*_t = \|X\|^*_{0,t}$ and $|X|^*_t = |X|^*_{0,t}$. The notation $X(t)$ and $X_t$ are used interchangeably. For a locally integrable function $f : \mathbb{R}_+ \to \mathbb{R}$ denote $\mathscr{I} f = \int_0^\cdot f(s)\,ds$. In case that $f$ is vector- or matrix-valued, $\mathscr{I} f$ is understood elementwise. For $y \in \mathbb{R}^k_+$ denote by $[\![y]\!]$ the element $y' \in \mathbb{Z}^k_+$ having $y'_i = \lfloor y_i \rfloor$ for $i = 1, \dots, k-1$, and $y'_k = \lfloor y_k \rfloor + \sum_{i=1}^k (y_i - \lfloor y_i \rfloor)$. Note that $e \cdot [\![y]\!] = e \cdot y$, and

$$\|y - [\![y]\!]\| \le 2k. \tag{1}$$

Finally, we need the following elementary graph-theoretic terms: A *tree* is an undirected, acyclic, connected graph; and the *diameter* of a graph is the distance between the two vertices that are farthest apart.

## 2. Setting and results.

2.1. *Queueing model.* The model under study has $I$ customer classes and $J$ service stations [see Figure 1(a)]. At each service station there are many independent servers of the same type. Each customer requires service only once and can be served indifferently by any server at the same station, but possibly at different rates at different stations. Only some stations can offer service to each class. When referring to the physical location of customers we say that they are in the buffer, or in the queue if they are not being served, and we say that they are in a certain station if they are served by a server of the corresponding type. There is one buffer per customer class and one station per server type. While a customer is in the queue it may abandon the system at rate depending on the class.

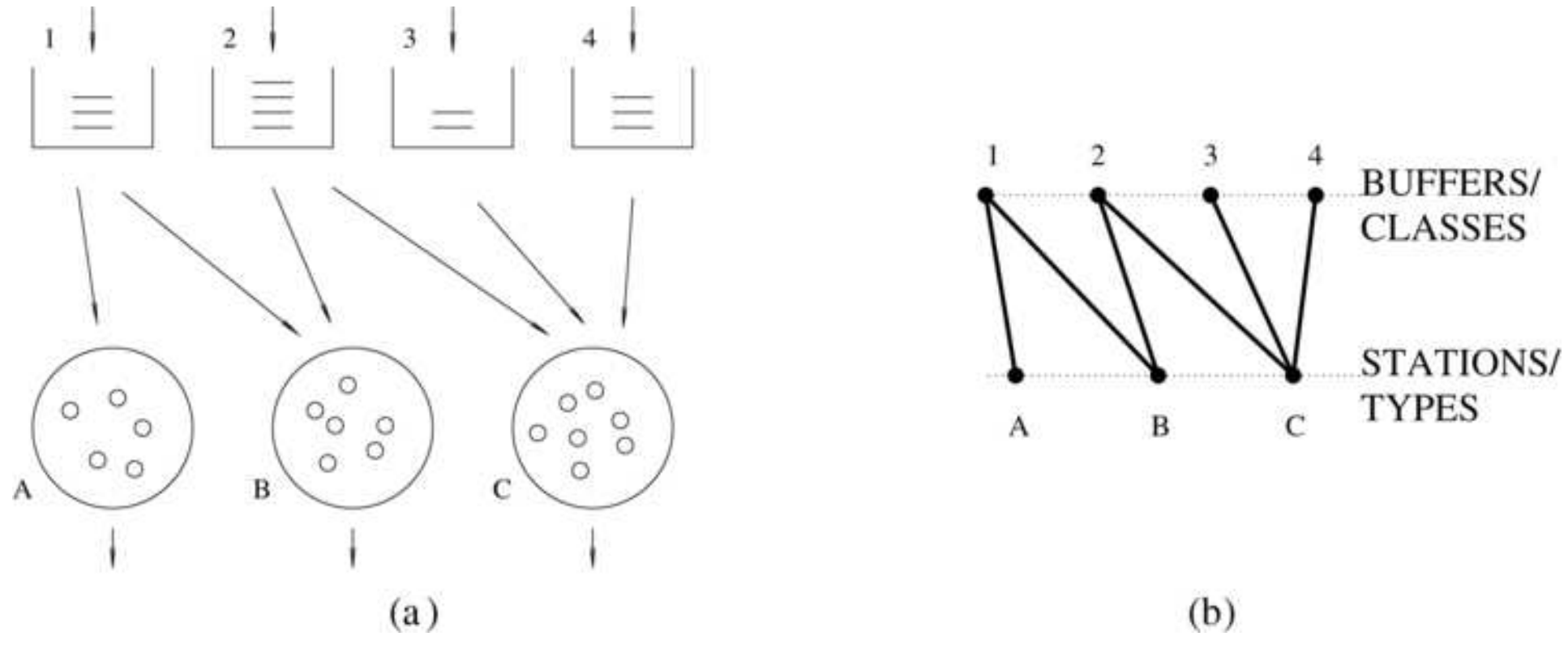

FIG. 1. (a) *A system with four classes and three stations.* (b) *The corresponding graph* $\mathcal{T}$.

A complete probability space $(\Omega, F, P)$ is given, supporting all stochastic processes defined below. Expectation with respect to $P$ is denoted by $E$. The processes will all be indexed by the parameter $n \in \mathbb{N}$, that will eventually be taken to infinity. The classes (buffers) are labelled as $1, \dots, I$ and the types (stations) as $I+1, \dots, I+J$:

$$\mathcal{I} = \{1, \dots, I\},$$
$$\mathcal{J} = \{I+1, \dots, I+J\}.$$

For $j \in \mathcal{J}$ let $N_j^n$ be the number of servers at station $j$. Let $X_i^n(t)$ denote the total number of class-$i$ customers in the system at time $t$. Let $Y_i^n(t)$ denote the number of class-$i$ customers in the queue at time $t$. Let $Z_j^n(t)$ denote the number of idle servers in station $j$ at time $t$. And let $\Psi_{ij}^n(t)$ denote the number of class-$i$ customers in station $j$ at time $t$. Let $X^n = (X_i^n)_{i\in\mathcal{I}}$, $Y^n = (Y_i^n)_{i\in\mathcal{I}}$, $Z^n = (Z_j^n)_{j\in\mathcal{J}}$, $\Psi^n = (\Psi_{ij}^n)_{i\in\mathcal{I}, i\in\mathcal{J}}$. Straightforward relations are expressed by the following equations, holding for all $i \in \mathcal{I}$, $j \in \mathcal{J}$ and $t \geq 0$:

$$Y_i^n(t) + \sum_{j\in\mathcal{J}} \Psi_{ij}^n(t) = X_i^n(t), \tag{2}$$

$$Z_j^n(t) + \sum_{i\in\mathcal{I}} \Psi_{ij}^n(t) = N_j^n, \tag{3}$$

$$Y_i^n(t) \geq 0, \qquad Z_j^n(t) \geq 0, \tag{4}$$

$$X_i^n(t) \geq 0, \qquad \Psi_{ij}^n(t) \geq 0. \tag{5}$$

To define arrival processes, let, for each $i \in \mathcal{I}$, $\{\check{U}_i(k), k \in \mathbb{N}\}$ be a sequence of strictly positive i.i.d. random variables with $E\check{U}_i(1) = 1$ and squared coefficient of variation $(E\check{U}_i(1))^{-2} \operatorname{Var}(\check{U}_i(1)) = C_{U,i}^2 \in [0, \infty)$. Assume also that the sequences are independent. Let

$$U_i^n(k) = \frac{1}{\lambda_i^n} \check{U}_i(k), \tag{6}$$

where $\lambda_i^n > 0$. With $\sum_1^0 = 0$, define

$$A_i^n(t) = \sup\left\{ l \geq 0 : \sum_{k=1}^{l} U_i^n(k) \leq t \right\}, \qquad t \geq 0. \tag{7}$$

It is assumed that the number of arrivals up to time $t$ is $A_i^n(t)$. Note that the first class-$i$ customer arrives at $U_i^n(1)$, and the time between the $(k-1)$st and $k$th arrival of class-$i$ customers is $U_i^n(k)$.

To model service times as exponential independent random variables, let $S^n_{ij}$, $i\in\mathcal{I}, j\in\mathcal{J}$, be Poisson processes with rate $\mu^n_{ij}\in[0,\infty)$ (where a zero-rate Poisson process is the zero process). These processes are assumed to be mutually independent, and independent of the arrival processes. Let $T^n_{ij}(t)$ denote the time up to $t$ devoted to a class-$i$ customer by a server, summed over all type-$j$ servers and note that

$$T^n_{ij}(t)=\int_0^t \Psi^n_{ij}(s)\,ds, \qquad i\in\mathcal{I}, j\in\mathcal{J}, t\ge 0.$$

The number of service completions of class-$i$ customers by all type-$j$ servers by time $t$ is $S^n_{ij}(T^n_{ij}(t))$. Note that this indeed models independent servers with exponential service times: Conditioned on $\Psi^n_{ij}=k$, the time until the next service, when $k$ servers are busy, is exponential with parameter proportional to $k$. Similarly, let $R^n_i(t)$ be Poisson processes of rate $\theta^n_i\in[0,\infty)$ and let $\mathring{T}^n_i(t)$ denote the time up to $t$ a class-$i$ customer spends in the queue, summed over all customers. Then

$$\mathring{T}^n_i(t)=\int_0^t Y^n_i(s)\,ds,$$

and we assume that $R^n_i(\mathring{T}^n_i(t))$ class-$i$ customers have abandoned until $t$. Throughout, the initial conditions, that we denote by $X^{0,n}_i:=X^n_i(0)$, are assumed to be deterministic. We have

$$X^n_i(t)=X^{0,n}_i+A^n_i(t)-\sum_j S^n_{ij}\biggl(\int_0^t \Psi^n_{ij}(s)\,ds\biggr)-R^n_i\biggl(\int_0^t Y^n_i(s)\,ds\biggr), \qquad i\in\mathcal{I}, t\ge 0. \tag{8}$$

The processes $A^n$, $S^n$ and $R^n$ will collectively be referred to as the primitive processes.

2.2. *Scheduling.* Scheduling decisions are made by continuously selecting $\Psi^n$, subject to appropriate constraints. Scheduling is regarded as *preemptive* if service to a customer can be stopped and resumed at a later time, possibly in a different station. Formally this is expressed by stating that the process $\Psi$ may be selected subject only to (2)–(12) holding. Note that according to this definition, customers can be moved instantaneously not only between a service station and the buffer, but also between different service stations that offer service to the corresponding class. Scheduling is regarded as *nonpreemptive* if every customer completes service with the server it is first assigned. More precisely, consider the processes $B^n_{ij}(t)$, $i\in\mathcal{I}, j\in\mathcal{J}$, where $B^n_{ij}(0)=0$, and $B^n_{ij}$ increases by $k$ each time $k$ class-$i$ jobs are moved to station $j$ from the buffer or from another station (to start or resume service), and decreases by $k$ each time $k$ such jobs are moved from station $j$ back to the buffer or to another station. Then $B^n_{ij}$ can be expressed as

$$B^n_{ij}(t) = \Psi^n_{ij}(t) - \Psi^{0,n}_{ij} + S^n_{ij}\left(\int_0^t \Psi^n_{ij}(s)\,ds\right). \tag{9}$$

To define nonpreemptive scheduling in terms of the model equations (2)–(12), we will require that $\Psi$ is selected subject to (2)–(12) and such that $B^n_{ij}$ are nondecreasing processes. This is summarized in the following.

DEFINITION 1. Let initial data $X^{0,n}$ and primitive processes $A^n$, $S^n$ and $R^n$ be given.

(i) We say that $\Psi^n$ is a *preemptive resume scheduling control policy* (P-SCP) if $\Psi_{ij}$ have cadlag paths ($i \in \mathcal{I}, j \in \mathcal{J}$), and there exist processes $X^n, Y^n$ and $Z^n$ such that (2)–(12) are met. $X^n$ is said to be the controlled process associated with initial data $X^{0,n}$ and P-SCP $\Psi^n$.

(ii) We say that $\Psi^n$ is a *nonpreemptive scheduling control policy* (N-SCP) if it is a P-SCP and $B^n_{ij}$ of (9) have nondecreasing paths.

We collectively refer to P-SCPs and N-SCPs as *scheduling control policies* (SCPs).

A word on terminology: The reader may notice that these definitions make nonpreemptive policies a special case of preemptive ones, unlike the usual use of these terms in the literature, where they are mutually exclusive. Indeed, the term "preemptive" is usually used to indicate that service to a customer is *necessarily* interrupted if a customer of higher priority arrives, whereas "nonpreemptive" is used when interruptions are not allowed. In the current context, however, the term "preemptive" indicates only that it is *possible* to interrupt service to a customer, and so it is natural that a policy under which no interruption occurs is regarded as a special case.

Let $\zeta = (X^{0,n}; n \in \mathbb{N})$ and $p = (\Psi^n, n \in \mathbb{N})$ denote a sequence of initial conditions, and respectively, SCPs. We denote by $P^p_\zeta$ the measure under which, for each $n$, $X^n$ is the controlled process associated with $X^{0,n}$ and $\Psi^n$. Expectation under $P^p_\zeta$ is denoted by $E^p_\zeta$.

We need a notion of SCPs that do not anticipate the future. Unlike in a Markovian setting, in the presence of renewal processes such a notion has to take into account that the time of the next arrival is correlated with information from the past, and hence should not be regarded as innovative information (see also [3]). Denote

$$\tau^n_i(t) = \inf\{u \ge t : A^n_i(u) - A^n_i(u-) > 0\}, \qquad i \in \mathcal{I}.$$

Set

$$\begin{aligned}\mathcal{F}_t^n = \sigma\{&A_i^n(s), S_{ij}^n(T_{ij}^n(s)), R_i^n(\mathring{T}_i^n(s)), \Psi_{ij}^n(s),\\ &X_i^n(s), Y_i(s), Z_j(s) : i \in \mathcal{I}, j \in \mathcal{J}, s \le t\}\end{aligned} \tag{10}$$

and

$$\begin{aligned}\mathcal{G}_t^n = \sigma\{&A_i^n(\tau_i^n(t) + u) - A_i^n(\tau_i^n(t)),\\ &S_{ij}^n(T_{ij}^n(t) + u) - S_{ij}^n(T_{ij}^n(t)),\\ &R_i^n(\mathring{T}_i^n(t) + u) - R_i^n(\mathring{T}_i^n(t)) : i \in \mathcal{I}, j \in \mathcal{J}, u \ge 0\}.\end{aligned} \tag{11}$$

DEFINITION 2. We say that a scheduling control policy is *admissible* if:

(i) for each $t$, $\mathcal{F}_t^n$ is independent of $\mathcal{G}_t^n$;
(ii) for each $i, j$ and $t$, the process $S_{ij}^n(T_{ij}^n(t) + \cdot) - S_{ij}^n(T_{ij}^n(t))$ [resp., $R_i^n(\mathring{T}_i^n(t) + \cdot) - R_i^n(\mathring{T}_i^n(t))$] is equal in law to $S_{ij}^n(\cdot)$ [$R_i^n(\cdot)$].

2.3. *Fluid and diffusion scaling.* We now introduce some assumptions regarding the asymptotic behavior of the parameters as $n \to \infty$. To this end, consider the graph $\mathcal{T}$ having vertex set $\mathcal{I} \cup \mathcal{J}$ with a vertex for each class and a vertex for each type, and an edge set $\mathcal{E}$, with an edge joining a class and a type if the corresponding service rate is nonzero:

$$\mathcal{E} = \{(i, j) \in \mathcal{I} \times \mathcal{J} : \mu_{ij}^n > 0\}.$$

It is assumed that the graph does not depend on the parameter $n$. We denote $i \sim j$ and $j \sim i$ if $(i, j) \in \mathcal{E}$. By assumption we have

$$\Psi_{ij}^n = 0, \qquad i \not\sim j. \tag{12}$$

A buffer-station pair $(i, j)$ is said to be an *activity* if $(i, j) \in \mathcal{E}$.

*Scaling of parameters.* There are constants $\lambda_i > 0$, $\theta_i \ge 0$, $i \in \mathcal{I}$, and $\mu_{ij} \ge 0$, $(i, j) \in \mathcal{I} \times \mathcal{J}$ where $\mu_{ij} > 0$ whenever $i \sim j$, such that

$$n^{-1}\lambda_i^n \to \lambda_i, \qquad \mu_{ij}^n \to \mu_{ij}, \qquad \theta_i^n \to \theta_i \tag{13}$$

and

$$n^{-1}N_j^n \to \nu_j.$$

Note that $\mu_{ij} = 0$ for $i \not\sim j$, since we had $\mu_{ij}^n = 0$ for $i \not\sim j$. Setting

$$\bar{\mu}_{ij} = \nu_j \mu_{ij}, \qquad i \in \mathcal{I}, j \in \mathcal{J},$$

a central assumption on the limit parameters, indicating that the sequence of systems is asymptotically critically loaded, is stated below (cf. [9, 10]).

*Linear program.* Minimize $\rho$ subject to

$$\sum_{j\in\mathcal{J}} \bar{\mu}_{ij}\xi_{ij} = \lambda_i, \qquad i \in \mathcal{I}, \tag{14}$$

$$\sum_{i\in\mathcal{I}} \xi_{ij} \le \rho, \qquad j \in \mathcal{J}, \tag{15}$$

$$\xi_{ij} \ge 0, \qquad i \in \mathcal{I}, j \in \mathcal{J}. \tag{16}$$

*Heavy traffic condition.* There exists a unique optimal solution $(\xi^*, \rho^*)$ to the linear program. Moreover, $\sum_{i\in\mathcal{I}} \xi^*_{ij} = 1$ for all $j \in \mathcal{J}$ (and, consequently, $\rho^* = 1$).

In the rest of this paper, $\xi^*_{ij}$ denotes the quantities from the above condition, and $x^* = (x^*_i)$, $\psi^* = (\psi^*_{ij})$, where

$$x^*_i = \sum_j \xi^*_{ij}\nu_j, \qquad \psi^*_{ij} = \xi^*_{ij}\nu_j. \tag{17}$$

We refer to the quantities $\xi^*_{ij}$, $x^*_i$ and $\psi^*_{ij}$ as *the static fluid model.*

In a fluid model where class-$i$ customers are served at station $j$ at rate $\bar{\mu}_{ij}$, the quantities $\sum_j \xi^*_{ij}\bar{\mu}_{ij}$ represent the actual service rate of class-$i$ customers when each station $j$ allocates a fraction $\xi^*_{ij}$ of its servers to class-$i$ jobs. The quantities $x^*_i$ represent the total mass of class-$i$ customers in the system. The heavy traffic condition indeed asserts that the system is critically loaded in the following sense. On one hand, (14) states that customers arriving at rates $\lambda_i$ can be processed by the system. On the other hand, under this condition it is impossible to have a $\tilde{\xi}$ satisfying (15) (with $\rho \le 1$) and (16), and a $\tilde{\lambda} \ne \lambda$, $\tilde{\lambda}_i \ge \lambda_i$, $i \in \mathcal{I}$, such that (14) holds for $\bar{\mu}, \tilde{\lambda}$ and $\tilde{\xi}$ (since this would violate the uniqueness of $\xi^*$). Thus, the system cannot process arrivals at rates greater than $\lambda_i$.

Following terminology from [10], an activity $(i, j) \in \mathcal{E}$ is said to be *basic* if $\xi^*_{ij} > 0$. This indicates that under the allocation matrix $\xi^*$, class-$i$ arrivals are actually processed by station $j$ in the fluid model. In this paper we assume that, in fact, all activities are basic. More precisely, we will assume:

Assumption 1. The heavy traffic condition holds, the graph $\mathcal{T}$ is connected and all activities are basic [i.e., $\xi^*_{ij} > 0$ for all $(i, j) \in \mathcal{E}$].

By a result of Williams ([14], Theorem 5.3 and Corollary 5.4), Assumption 1 implies that

$$\text{the graph } \mathcal{T} \text{ is a tree.} \tag{18}$$

As a result, Assumption 1 can equivalently be stated as: *The heavy traffic condition holds, the graph $\mathcal{T}$ is a tree and all activities are basic.* The assumption that $\mathcal{T}$ is a tree was exploited in the PDE analysis in [2]. We have more to say about Assumption 1 in Section 2.7.

Second-order scaling assumptions are as follows.

*Scaling of parameters, continued.* There are constants $\hat\lambda_i$, $\hat\mu_{ij}\in\mathbb{R}$, $i\in\mathcal{I}$, $j\in\mathcal{J}$, such that

$$\hat\lambda_i^n := n^{1/2}(n^{-1}\lambda_i^n-\lambda_i)\to\hat\lambda_i, \qquad \hat\mu_{ij}^n := n^{1/2}(\mu_{ij}^n-\mu_{ij})\to\hat\mu_{ij}, \tag{19}$$

$$n^{1/2}(n^{-1}N_j^n-\nu_j)\to 0. \tag{20}$$

*Scaling of initial conditions.* There are constants $x_i$, $y_i$, $z_j$, $\psi_{ij}$ satisfying $y_i+\sum_j\psi_{ij}=x_i$, $z_j+\sum_i\psi_{ij}=0$, $y_i\ge 0$, and $z_j\ge 0$, $i\in\mathcal{I}$, $j\in\mathcal{J}$, such that the deterministic initial conditions satisfy

$$\begin{aligned}\hat X_i^{0,n} &:= n^{-1/2}(X_i^{0,n}-nx_i^*)\to x_i,\\ \hat Y_i^{0,n} &:= n^{-1/2}Y_i^{0,n}\to y_i,\end{aligned} \tag{21}$$

$$\begin{aligned}\hat Z_j^{0,n} &:= n^{-1/2}Z_j^{0,n}\to z_j,\\ \hat\Psi_{ij}^{0,n} &:= n^{-1/2}(\Psi_{ij}^{0,n}-\psi_{ij}^*n)\to\psi_{ij}.\end{aligned} \tag{22}$$

The processes rescaled at the fluid level are defined as

$$\bar X_i^n(t)=n^{-1}X_i^n(t), \qquad \bar Y_i^n(t)=n^{-1}Y_i^n(t),$$
$$\bar Z_j^n(t)=n^{-1}Z_j^n(t), \qquad \bar\Psi_{ij}^n(t)=n^{-1}\Psi_i^n(t).$$

The primitive processes are centered about their means and rescaled at the diffusion level as

$$\begin{aligned}\hat A_i^n(t) &= n^{-1/2}(A_i^n(t)-\lambda_i^n t),\\ \hat S_{ij}^n(t) &= n^{-1/2}(S_{ij}^n(nt)-n\mu_{ij}^n t),\\ \hat R_i^n(t) &= n^{-1/2}(R_i^n(nt)-n\theta_i^n t).\end{aligned} \tag{23}$$

Similarly, the state processes are centered about the static fluid model and rescaled:

$$\hat X_i^n(t)=n^{-1/2}(X_i^n(t)-nx_i^*), \tag{24}$$

$$\begin{aligned}\hat Y_i^n(t) &= n^{-1/2}Y_i^n(t),\\ \hat Z_j^n(t) &= n^{-1/2}Z_j^n(t),\end{aligned} \tag{25}$$

$$\hat\Psi_{ij}^n(t)=n^{-1/2}(\Psi_{ij}^n(t)-\psi_{ij}^*n). \tag{26}$$

The relations (2), (3) and (4) take the new form

$$\hat Y_i^n+\sum_j\hat\Psi_{ij}^n=\hat X_i^n, \qquad i\in\mathcal{I}, \tag{27}$$

$$\hat Z_j^n+\sum_i\hat\Psi_{ij}^n=0, \qquad j\in\mathcal{J}, \tag{28}$$

$$\hat Y_i^n,\hat Z_j^n\ge 0, \qquad i\in\mathcal{I}, j\in\mathcal{J}. \tag{29}$$

Using the definitions above one finds that (8) can be written as follows (for the reader's convenience this development is provided in the Appendix):

$$\hat{X}_i^n(t) = \hat{X}_i^{0,n} + r_i\hat{W}_i^n(t) + \ell_i^n t - \sum_j \mu_{ij}^n \int_0^t \hat{\Psi}_{ij}^n(s)\,ds - \theta_i^n \int_0^t \hat{Y}_i^n(s)\,ds, \tag{30}$$

where

$$\begin{aligned} r_i\hat{W}_i^n(t) &= \hat{A}_i^n(t) - \sum_j \hat{S}_{ij}^n\left(\int_0^t \bar{\Psi}_{ij}^n(s)\,ds\right) - \hat{R}_i^n\left(\int_0^t \bar{Y}_i^n(s)\,ds\right), \\ \ell_i^n &= \hat{\lambda}_i^n - \sum_j \hat{\mu}_{ij}^n \psi_{ij}^*. \end{aligned} \tag{31}$$

With (19) we have

$$\lim_n \ell_i^n = \ell_i := \hat{\lambda}_i - \sum_j \hat{\mu}_{ij}\psi_{ij}^*. \tag{32}$$

One is free to choose the values of $r_i$, and it will be convenient to choose them so that, with the formal substitution $\bar{\Psi}_{i,j}^n = \psi_{i,j}^*$, $\bar{Y}_i^n = 0$, one has

$$\lim_n E[(\hat{W}_i^n(1))^2] = 1. \tag{33}$$

Namely, $r_i = (\lambda_i C_{U,i}^2 + \lambda_i)^{1/2}$. Denote also $\ell = (\ell_1, \ldots, \ell_I)'$, $r = \mathrm{diag}(r_i)$.

The results of this paper are concerned with constructing sequences of SCPs that, in an appropriate sense, minimize the limit as $n \to \infty$ of cost of the following form:

$$E\int_0^\infty e^{-\gamma t}\tilde{L}(\hat{X}_t^n, \hat{\Psi}_t^n)\,dt. \tag{34}$$

2.4. *Joint work conservation.* A policy is said to be *work conserving* if it does not allow for a server to idle while a customer that it can serve is in the queue. We can express this condition as

$$Y_i^n(t) \wedge Z_j^n(t) = 0 \qquad \forall\, i \sim j, t \ge 0.$$

In the current context one can consider a stronger condition for preemptive policies. Recall that if preemption is allowed, customers of each class can be moved between the queue and the various stations that offer service to them. In a preemptive policy, controls $\Psi^n(t)$ that can be applied at time $t$ correspond to different rearrangements of the customers $X^n(t)$ in the stations and buffers. Let $\mathcal{X}^n$ denote the set of all possible values of $X^n(t)$ for which there is a rearrangement of customers with the property: *either there are no customers in the queue, or no server in the system is idle.* This property is expressed as

$$e \cdot Y^n(t) \wedge e \cdot Z^n(t) = 0. \tag{35}$$

We shall say that a preemptive policy is *jointly work conserving* if it is work conserving and, in addition, for every $t$, if $X^n(t) \in \mathcal{X}^n$, then customers are arranged according to (35). Some motivation is provided in Section 2.7.

Let $\alpha_0 > 0$ denote the constant

$$\alpha_0 = (4C_G)^{-1} \min_{i,j\,:\,i\sim j} \psi^*_{ij}. \tag{36}$$

It will be shown (see Lemma 3) that

$$\|X^n(t) - nx^*\| \le \alpha_0 n \quad \text{implies} \quad X^n(t) \in \mathcal{X}^n \tag{37}$$

(note that the value of $t$ is irrelevant in the discussion). This fact is significant for the following reason. Under appropriate conditions, it will be established that the processes $\hat{X}^n$ are tight. As a result, the condition on the left-hand side of (37) holds for every $t \in [0,T]$ (for arbitrary $T$), with probability approaching 1 as $n \to \infty$. Thus it is nearly always the case that a rearrangement of customers according to (35) is possible. These details will be justified in the process of proving our result.

2.5. *The diffusion control problem.* We take limits as $n \to \infty$ in (27), (28), (30) and (35). The process $\hat{W}^n$ converges to a standard Brownian motion, and denoting the weak limit of $(\hat{X}^n, \hat{Y}^n, \hat{Z}^n, \hat{W}^n, \hat{\Psi}^n)$ by $(X, Y, Z, W, \Psi)$, we obtain the equation below (at this point this is meant as a formal step only; however, see Proposition 3):

$$X_i(t) = x_i + \tilde{W}_i(t) - \sum_j \mu_{ij} \int_0^t \Psi_{ij}(s)\,ds - \theta_i \int_0^t Y_i(s)\,ds, \qquad i \in \mathcal{I}, \tag{38}$$

where $\tilde{W}_i(t) = r_i W_i(t) + \ell_i t$, $W$ is a standard Brownian motion and

$$\sum_j \Psi_{ij} = X_i - Y_i, \qquad i \in \mathcal{I}, \tag{39}$$

$$\sum_i \Psi_{ij} = -Z_j, \qquad j \in \mathcal{J}, \tag{40}$$

$$e \cdot Y \wedge e \cdot Z = 0. \tag{41}$$

Relations (38)–(41) above can be written in the convenient form

$$dX = b(X, U)\,dt + r\,dW.$$

To this end, note that by (39) and (40), $e \cdot X = e \cdot Y - e \cdot Z$, and thus by (41)

$$e \cdot Y = (e \cdot X)^+, \qquad e \cdot Z = (e \cdot X)^-.$$

Hence $Y$ and $Z$ can be represented in terms of the process $e\cdot X$ and an additional process $U$ as

$$Y_i(t)=(e\cdot X(t))^+u_i(t),\qquad Z_j(t)=(e\cdot X(t))^-v_j(t),\qquad i\in\mathcal{I},j\in\mathcal{J},\tag{42}$$

where $U(t)=(u(t),v(t))$ takes values in

$$\mathbb{U}:=\{(u,v)\in\mathbb{R}^{I+J}:u_i,v_j\ge 0,i\in\mathcal{I},j\in\mathcal{J},\ e\cdot u=e\cdot v=1\}.$$

The following is shown in Proposition 7 of [2].

Lemma 1. *Let Assumption 1 hold. Then given $\alpha_i,\beta_j\in\mathbb{R}$, $i\in\mathcal{I},j\in\mathcal{J}$ satisfying $e\cdot\alpha=e\cdot\beta$, there exists a unique solution $\psi=(\psi_{ij},i\in\mathcal{I},j\in\mathcal{J})\in\mathbb{R}^{IJ}$ to the set of equations*

$$\begin{aligned}&\sum_j\psi_{ij}=\alpha_i,\qquad i\in\mathcal{I},\sum_i\psi_{ij}=\beta_j,\qquad j\in\mathcal{J},\\&\psi_{ij}=0\qquad\text{for } i\not\sim j.\end{aligned}\tag{43}$$

As a result there is a map, denoted throughout by $G:D_G\to\mathbb{R}^{IJ}$,

$$D_G:=\{(\alpha,\beta)\in\mathbb{R}^{I+J}:e\cdot\alpha=e\cdot\beta\},$$

such that $\psi=G(\alpha,\beta)$. Let also

$$C_G=\sup\Big\{\max_{ij}|G_{ij}(\alpha,\beta)|:(\alpha,\beta)\in D_G,\|\alpha\|\vee\|\beta\|\le 1\Big\}.\tag{44}$$

Clearly the map is linear on $D_G$. Applying Lemma 1 to (39), (40), and by (42) it follows that

$$\Psi=G(X-Y,-Z)=G(X-(e\cdot X)^+u,-(e\cdot X)^-v)=:\hat G(X,U).\tag{45}$$

Letting

$$b_i(X,U)=-\sum_{j\in\mathcal{J}}\mu_{ij}\hat G_{ij}(X,U)-\theta_i(e\cdot X)^+u_i+\ell_i,\tag{46}$$

we can now write (38) as

$$X(t)=x+rW(t)+\int_0^t b(X(s),U(s))\,ds,\qquad 0\le t<\infty.\tag{47}$$

We summarize the equivalence between the two representations in what follows.

Lemma 2. *If (38)–(41) hold, then (47) holds with some $\mathbb{U}$-valued $U$. Conversely, if (47) holds (with $U$ taking values in $\mathbb{U}$), then one can find $\Psi,Y,Z$ such that (38)–(41) hold. In both cases, $\tilde W(t)=rW(t)+\ell t$.*

PROOF. We have already proved the first statement of the result. To see the converse, write $U=(u,v)$, let $Y=(e\cdot X)^+u$, $Z=(e\cdot X)^-v$, and $\Psi=G(X-Y,-Z)$. Equations (39), (40) automatically hold, and (45) implies (38). $\square$

DEFINITION 3. We call $\pi=(\Omega,F,(F_t),P,U,W)$ an *admissible system* if $(\Omega,F,(F_t),P)$ is a complete filtered probability space, $U$ is a $\mathbb{U}$-valued, $(F_t)$-progressively measurable process, and $W$ is a standard $I$-dimensional $(F_t)$-Brownian motion. The process $U$ is said to be a *control* associated with $\pi$. $X$ is said to be a *controlled process* associated with initial data $x\in\mathbb{R}^I$ and an admissible system $\pi$, if it is a continuous sample paths, $(F_t)$-adapted process such that $\int_0^t|b(X(s),U(s))|\,ds<\infty$, $t\ge 0$, $P$-a.s., and (47) holds $P$-a.s.

For any $x\in\mathbb{R}^I$ and any admissible system $\pi$ there exists a controlled process $X$, unique in the strong sense (cf. [2]). With an abuse of notation we sometimes denote the dependence on $x$ and $\pi$ by writing $P_x^\pi$ in place of $P$ and $E_x^\pi$ in place of $E$. We denote by $\Pi$ the class of all admissible systems.

Given a constant $\gamma>0$ and a function $\tilde{L}$, the cost of interest for the queueing system is given by (34). It is convenient to perform change of variables from $(X,\Psi)$ to $(X,U)$. To this end, define $L$ as

$$L(X,U)=\tilde{L}(X,G(X-(e\cdot X)^+u,-(e\cdot X)^-v)),\qquad X\in\mathbb{R}^I,U=(u,v)\in\mathbb{U}. \tag{48}$$

Our conditions on $\tilde{L}$ (given mostly via conditions on $L$) are stated below. They were required to carry out the PDE analysis in [2].

ASSUMPTION 2. (i) $L(x,U)\ge 0$, $(x,U)\in\mathbb{R}^I\times\mathbb{U}$.

(ii) The mapping $(x,\psi)\mapsto\tilde{L}(x,\psi)$ is continuous. In particular, the mapping $(x,U)\mapsto L(x,U)$ is continuous.

(iii) There is $\varrho\in(0,1)$ such that for any compact $A\subset\mathbb{R}^I$,

$$|L(x,U)-L(y,U)|\le c\|x-y\|^\varrho$$

holds for $U\in\mathbb{U}$ and $x,y\in A$, where $c$ depends only on $A$.

(iv) There are constants $c>0$ and $m_L\ge 1$ such that $L(x,U)\le c(1+\|x\|^{m_L})$, $U\in\mathbb{U}$, $x\in\mathbb{R}^I$.

Consider the cost

$$C(x,\pi)=E_x^\pi\int_0^\infty e^{-\gamma t}L(X(t),U(t))\,dt,\qquad x\in\mathbb{R}^I,\pi\in\Pi.$$

Define the value function as

$$V(x) = \inf_{\pi\in\Pi} C(x,\pi).$$

The HJB equation for the problem is

$$\mathcal{L}f + H(x, Df) - \gamma f = 0, \tag{49}$$

where $D$ denotes the gradient, $\mathcal{L} = (1/2)\sum_i r_i^2\, \partial^2/\partial x_i^2$, and

$$H(x,p) = \inf_{U\in\mathbb{U}} [b(x,U)\cdot p + L(x,U)]. \tag{50}$$

The equation is considered on $\mathbb{R}^I$ with the growth condition

$$\exists\, C, m \qquad |f(x)| \le C(1 + \|x\|^m), \qquad x \in \mathbb{R}^I. \tag{51}$$

DEFINITION 4. Let $x\in\mathbb{R}^I$ be given. We say that a measurable function $h:\mathbb{R}^I\to\mathbb{U}$ is a *Markov control policy* if there is an admissible system $\pi$ and a controlled process $X$ corresponding to $x$ and $\pi$, such that $U_s = h(X_s)$, $s\ge 0$, $P$-a.s. We say that an admissible system $\pi$ is *optimal* for $x$, if $V(x) = C(x,\pi)$. We say that a Markov control policy is optimal for $x$ if at least one of the admissible systems corresponding to it is optimal.

We note in passing that the last part of the above definition is equivalent to optimality of *all* corresponding admissible systems, due to weak uniqueness of solutions to $X(t) = x + rW(t) + \int_0^t b(X(s), h(X(s)))\,ds$ (cf. Proposition 5.3.10 of [11]).

The assumption below is needed for large time estimates on the state processes. Unless $L$ is bounded [as in part (iii)], such estimates are required for both the DC problem and the SC problem, and they are only known to hold under parts (i) and (ii) below. This is explained in the remark in Section 3.3.

ASSUMPTION 3. Either (i), (ii) or (iii) below holds.

(i) For $(i,j)\in\mathcal{E}$, $\mu_{ij}$ depends only on $i$; or for $(i,j)\in\mathcal{E}$, $\mu_{ij}$ depends only on $j$. In addition, $\theta_i = 0$, $i\in\mathcal{I}$.

(ii) The tree $\mathcal{T}$ is of diameter 3 at most. In addition, for all $(i,j)\in\mathcal{E}$, $\theta_i \le \mu_{ij}$.

(iii) The function $L$ is bounded.

The following is proved in [2].

THEOREM 1. *Let Assumptions 1, 2 and 3 hold. Then the value $V$ solves (49), (51), and there exists a Markov control policy $h^*:\mathbb{R}^I\to\mathbb{U}$ that is optimal for all $x\in\mathbb{R}^I$. In cases* (i) *and* (ii) [*resp.,* (iii)] *of Assumption 3, uniqueness holds in $C^2_{\mathrm{pol}}$* [*resp.,* $C^2_b$].

2.6. *Asymptotically optimal SCPs.* To state our main result we need to introduce SCPs defined via the function $h^*$ of Theorem 1. Write $h^* = (h_1^*, h_2^*)$ where $(u, v) = h^*(x) \Leftrightarrow u = h_1^*(x), v = h_2^*(x)$. Theorem 1 states that, for the DC problem, there is an optimal pair $(U, X)$ for which $U(t) = h^*(X(t))$, or equivalently $u(t) = h_1^*(X(t))$, $v(t) = h_2^*(X(t))$. In view of (42) and (45), this can be written as

$$Y(t) = (e \cdot X(t))^+ h_1^*(X(t)), \qquad Z(t) = (e \cdot X(t))^- h_2^*(X(s)),$$
$$\Psi(t) = G(X(t) - Y(t), -Z(t)).$$

Let us define

$$\check{Y}^n(t) = (e \cdot \hat{X}^n(t))^+ h_1^*(\hat{X}^n(t)), \qquad \check{Z}^n(t) = (e \cdot \hat{X}^n(t))^- h_2^*(\hat{X}^n(t)), \tag{52}$$
$$\check{\Psi}^n(t) = G(\hat{X}^n(t) - \check{Y}^n(t), -\check{Z}^n(t)). \tag{53}$$

In analogy with the DC problem, one would like to define a P-SCP simply by setting

$$\hat{\Psi}^n(t) = \check{\Psi}^n(t). \tag{54}$$

However, this is impossible, because the components of the corresponding process $\Psi^n(t) = \psi^* n + n^{1/2} \hat{\Psi}^n(t)$ [cf. (26)] represent number of customers, and must be $\mathbb{Z}_+$-valued. Instead, a P-SCP will be defined in a way that the equality (54) will hold "roughly." To this end, consider two cases.

*Case* 1. If $\|X^n(t) - nx^*\| \le \alpha_0 n$, set

$$\begin{aligned} Y^n(t) &= [\![ n^{1/2} \check{Y}^n(t) ]\!], \\ Z^n(t) &= [\![ n^{1/2} \check{Z}^n(t) ]\!], \\ \Psi^n(t) &= G(X^n(t) - Y^n(t), N^n - Z^n(t)). \end{aligned} \tag{55}$$

*Case* 2. If $\|X^n(t) - nx^*\| > \alpha_0 n$, set $\Psi^n = F_n(X^n(t))$ where, for each $n$, $F_n$ is a fixed function chosen in such a way that the resulting process $\Psi^n$ is jointly work conserving (as defined in Section 2.4). Other than that the choice of $F_n$ is immaterial.

The choice of $\Psi^n(t)$ in Case 1 is in accordance with joint work conservation. Namely, (35) always holds on $\{\|X^n(t) - nx^*\| \le \alpha_0 n\}$. To see this note that by (52), $e \cdot \check{Y}^n(t) \wedge e \cdot \check{Z}^n(t) = 0$, hence by (55), $e \cdot Y^n(t) \wedge e \cdot Z^n(t) = 0$. Also, since by definition $X^n, Y^n$ and $Z^n$ are integer valued, so is $\Psi^n$ as follows from the proof of Proposition 1 of [2]. It remains to show that setting $\Psi^n(t)$ as in (55) meets relation (5) of the model, namely that $\Psi_{ij}^n(t) \ge 0$ for all $i, j$. The proof of this fact is deferred to Section 3.1 (cf. Lemma 3 and the remark that follows). Denote the resulting sequence of P-SCPs by $p^*$.

The proposed N-SCP is described next. The idea is still to keep $\hat{\Psi}^n$ close to $\check{\Psi}^n$, but it is no longer possible to assign the values of $\hat{\Psi}^n$ in a direct way: $\hat{\Psi}^n$ is indirectly affected by the job assignments. Heuristically, the policy we propose is a kind of tracking mechanism, in which $\hat{\Psi}^n$ tracks $\check{\Psi}^n$, as follows:

(a) Compute $\check{\Psi}^n(t)$ via (52), (53).
(b) Declare activities $(i,j) \in \mathcal{E}$ with $\hat{\Psi}^n_{ij}(t) > \check{\Psi}^n_{ij}(t)$ as "overpopulated."
(c) Block overpopulated activities. That is to say that no new job assignments are allowed on an activity as long as it is overpopulated. As a result, it is anticipated that the population in these activities drops rapidly, and the population in the other activities increases rapidly.
(d) Obtain $\hat{\Psi}^n(t) \sim \check{\Psi}^n(t)$.

We turn to the precise definition of the N-SCP. For a technical reason, we must change the definition (52) of $\check{Y}^n(t)$ and $\check{Z}^n(t)$ so that the function $h^*$ is replaced by another function for large values of $t$. Namely, let $h^0 = (h^0_1, h^0_2)$ where $h^0_1(x) = e_1$ and $h^0_2(x) = e_{I+1}$ for all $x$. Instead of (52), we let $\check{Y}^n(t)$ and $\check{Z}^n(t)$ be defined as

$$\check{Y}^n(t) = (e \cdot \hat{X}^n(t))^+ h_1(t, \hat{X}^n(t)), \qquad \check{Z}^n(t) = (e \cdot \hat{X}^n(t))^- h_2(t, \hat{X}^n(t)),$$

where

$$h(t,x) = (h_1(t,x), h_2(t,x)) = \begin{cases} h^*(x), & t < \Theta_n, \\ h^0(x), & t \geq \Theta_n, \end{cases} \tag{56}$$

and $\Theta_n$ is a random time to be defined below. The definition of $\check{\Psi}^n$ is still via (53).

Next, we would like the following to hold. Given an activity $(i,j)$ and a time interval $[s,t)$, if $\hat{\Psi}^n_{ij} > \check{\Psi}^n_{ij}$ holds on $[s,t)$, then no routings take place on the activity throughout this interval:

$$\hat{\Psi}^n_{ij} > \check{\Psi}^n_{ij} \text{ on } [s,t) \quad \text{implies} \quad B^n_{ij} \text{ is constant on } [s,t). \tag{57}$$

On the other hand, when there is a class-$i$ customer in the queue and there are stations $j \sim i$ with idle servers and such that $(i,j)$ is not blocked, the customer is instantaneously routed to one of these stations. If there are no such stations, the customer stays in the queue. It is not hard to see that this property can be expressed as follows. For every activity $(i,j)$,

$$\hat{\Psi}^n_{ij}(t) \leq \check{\Psi}^n_{ij}(t) \quad \text{implies} \quad Y^n_i(t) \wedge Z^n_j(t) = 0. \tag{58}$$

The selection of an activity among the nonblocked activities through which to route a customer does not turn out to be important. However, care must be taken when two (or more) customers are routed instantaneously. Therefore it remains to show that one can always perform instantaneous routings meeting (58). This is deferred to Section 3.1 (cf. Lemma 4). Finally, we define $\Theta_n$ as

$$\Theta_n = \inf\Big\{t : \max_{i\sim j}[\hat{\Psi}^n_{ij}(t_-) - F^*_{ij}(\hat{X}^n(t_-))] \geq b_0\Big\}, \tag{59}$$

where

$$F^*(x) = G(x - (e\cdot x)^+ h^*_1(x), -(e\cdot x)^- h^*_2(x)),$$

and $b_0$ denotes the deterministic, finite constant

$$b_0 = 2 + \max_{i\sim j}\sup_n(\hat{\Psi}^n_{ij}(0) - \check{\Psi}^n_{ij}(0))^+. \tag{60}$$

A sequence of N-SCPs satisfying the above properties will be denoted by $p'$. We will also consider a sequence of N-SCPs defined exactly as above, except that for each $n$, the function $h^*$ will be replaced by a function $h_n$, depending on $n$. Given a sequence $\{h_n\}$ of functions $h_n : \mathbb{R}^I \to \mathbb{U}$, $n \in \mathbb{N}$, we denote the corresponding sequence of N-SCPs as $p'(\{h_n\})$.

The estimates used to establish asymptotic optimality will require some assumptions on the interarrival times. The two parts of the assumption below correspond to different parts of the result.

Assumption 4. There is a constant $m_A$ such that $E(\check{U}_i(1))^{m_A} < \infty$, $i \in \mathcal{I}$, satisfying either of the following:

(i) $m_A > 2m_L$ (where $m_L$ is as in Assumption 2);
(ii) $m_A(m_A - 2)(5m_A - 2)^{-1} > m_L$.

For a sequence $\zeta$ of initial conditions and a sequence $p$ of SCPs, denote

$$\underline{V}(\zeta, p) = \liminf_{n\to\infty} E^p_\zeta \int_0^\infty e^{-\gamma t}\tilde{L}(\hat{X}^n_t, \hat{\Psi}^n_t)\,dt,$$

$$\overline{V}(\zeta, p) = \limsup_{n\to\infty} E^p_\zeta \int_0^\infty e^{-\gamma t}\tilde{L}(\hat{X}^n_t, \hat{\Psi}^n_t)\,dt.$$

Our main result is the following.

Theorem 2. *Let Assumptions 1, 2 and 3 hold. Let $\zeta$ be a sequence of initial conditions $(X^{0,n}; n \in \mathbb{N})$ such that $\hat{X}^{0,n} = n^{-1/2}(X^{0,n} - nx^*) \to x \in \mathbb{R}^I$. Then items* (i) *and* (ii) *below hold under Assumption 4*(i)*, and items* (iii) *and* (iv) *hold under Assumption 4*(ii)*.*

(i) *For any sequence $p$ of jointly work-conserving admissible P-SCPs, $\underline{V}(\zeta, p) \geq V(x)$.*

(ii) *The sequence $p^*$ of jointly work-conserving admissible P-SCPs satisfies $\overline{V}(\zeta, p^*) \leq V(x)$.*

(iii) *Provided that* $h^*$ *(of Theorem 1) is locally Hölder on* $\{\xi \in \mathbb{R}^I : e \cdot \xi \neq 0\}$*, the sequence* $p'$ *of admissible N-SCPs satisfies* $\overline{V}(\zeta, p') \leq V(x)$.

(iv) *Provided that the mapping* $U \mapsto L(\xi, U)$ *is convex on* $\mathbb{U}$ *for every* $\xi \in \mathbb{R}^I$*, there exists a sequence* $\{h_n\}$ *of functions mapping* $\mathbb{R}^I$ *to* $\mathbb{U}$ *such that the corresponding sequence* $p'' := p'(\{h_n\})$ *of admissible N-SCPs satisfies* $\overline{V}(\zeta, p'') \leq V(x)$.

The combination of items (i) and (ii) establishes the asymptotic optimality of the sequence $p^*$, as well as the validity of the DC problem as the correct asymptotic description of the problem under preemption. Items (i), (iii) and (iv) establish similar consequences for the nonpreemptive problem.

To demonstrate the horizon of the different parts of the result, we mention some possible cost functions associated with queue length, and note in passing that one can also consider, for example, costs associated with the number of idle servers, and, via a transformation developed in [3], abandonment count.

Recall that $\hat{Y}_i^n$ denote normalized queue lengths, and that $\hat{Y}_i^n = (e \cdot \hat{X}^n)^+ u_i^n$. For simplicity write these quantities as $Y_i$, $X_i$ and $u_i$. Then a cost of the form $\sum_i c_i Y_i^\alpha$, where $c_i, \alpha > 0$, can be written as

$$L(X, U) = ((e \cdot X)^+)^\alpha \sum_i c_i u_i^\alpha.$$

We see that Assumption 2 is met and therefore items (i) and (ii) of Theorem 2 are in force. For item (iii), it follows from Proposition 3 of [3] that the Hölder hypothesis on $h^*$ is met for the above cost with $\alpha \geq 2$. For item (iv), clearly $\alpha \geq 1$ suffices. Next, our treatment of the case where $L$ is bounded [cf. Assumption 3(iii)] allows to relax much of what is assumed otherwise [i.e., Assumption 3(i) or (ii)]. A truncated version of the above cost function, that is,

$$L(X, U) = \left[((e \cdot X)^+)^\alpha \sum_i c_i u_i^\alpha\right] \wedge M,$$

is an example in which items (i) and (ii) apply under Assumption 3(iii). As far as item (iv) is concerned, the convexity condition imposes the choice $\alpha = 1$ and $c = e$, which amounts to the truncated sum of the queue lengths $(e \cdot Y) \wedge M$, expressed as

$$L(X, U) = (e \cdot X)^+ \wedge M.$$

Clearly, any other positive, bounded Hölder function of $X$ meets the assumptions of item (iv).

As a final comment we mention that the PDE (49), on which the above proposed SCPs depend, can be solved numerically. The complexity of numerical schemes depend mainly on the dimension of the underlying Euclidean space. The PDE (49) is of dimension $I$, and to realize the dimensionality reduction obtained by our results as far as *nonpreemptive* SC problems are concerned, note that the state process for such a problem in the Markovian case should consist of $(Y^n, \Psi^n)$. Since the number of activities under the treelike assumption is equal to $I + J - 1$, the reduction is from dimension $2I + J - 1$ to $I$.

2.7. *Discussion.* We provide some motivation on two central assumptions that we have made. First, we present an example that demonstrates heuristically why, in presence of nonbasic activities, the limit behavior is expected to be different for preemptive and nonpreemptive scheduling. The example consists of the simplest system having a nonbasic activity. The observation and example are due to Mandelbaum and Reiman [12].

EXAMPLE 1. Consider a system with two customer classes and two server types as depicted in Figure 2. Servers of type $A$ can serve both classes and servers of type $B$ can only serve class-2 customers (in the examples we refer to stations as $A, B, \ldots$ rather than as numbers in $\mathcal{J}$). Consider first a case where all activities are basic. For example, let $\nu_1 = \nu_2 = 1$, $\mu_{1A} = \mu_{2A} = \mu_{2B} = 1$, and $\lambda_1 = 1/2$, $\lambda_2 = 3/2$. Then there is a unique optimal solution to the linear program, namely: $(\xi^*_{1A}, \xi^*_{2A}, \xi^*_{2B}) = (1/2, 1/2, 1)$. In a fluid model this means that half of the servers at station $A$ serve class 1 and all other servers serve class 2. In the queueing model with large $n$, under any policy that keeps the deviations from the fluid model at the desired level $O(\sqrt{n})$, the number of class-1 and the number of class-2 customers in station $A$ are $\Psi^n_{1A}, \Psi^n_{2A} = n/2 + O(\sqrt{n})$, while the number of class-2 customers in station $B$ is $\Psi^n_{2B} = n + O(\sqrt{n})$. In a system with preemption the policy dynamically selects the values of these processes $\Psi^n$. The control $\Psi$ that appears in the diffusion model shows up as a weak limit of the appropriately normalized, *centered* versions of the $\Psi^n$ process. It therefore reflects the $O(\sqrt{n})$ deviations referred to above. Now, in order for a system without preemption to be governed by the same diffusion model, it should be possible for $\Psi^n$ to track closely the diffusion control process $\Psi$. In particular, it should be possible for the population in any activity to change by $O(\sqrt{n})$ in $o(1)$ units of time (this is the heuristic behind the proofs carried out in Section 3). Recall that the population in this activity is about $n/2$. If routing to this activity is stopped, its population will decrease due to service completions by $O(\sqrt{n})$ in $o(1)$ units of time. On the other hand, if routing to this activity is maintained and routing to the other activities is stopped, the population at activity $(2,A)$ will build up by $O(\sqrt{n})$ in $o(1)$ units of time. This demonstrates how to rapidly decrease or increase the population in this activity, so as to allow $\Psi^n$ to track a desired control process. Consider now a case where $\mu$ and $\nu$ are as before, but $\lambda_1=\lambda_2=1$. The unique solution to the linear program is now $(\xi^*_{1A},\xi^*_{2A},\xi^*_{2B})=(1,0,1)$, and therefore activity $(2,A)$ is nonbasic. The queueing model will have only $O(\sqrt{n})$ customers in this activity. As a result, it is impossible to reduce the population in this activity in $o(1)$ units of time so as to maintain a "good" tracking mechanism. This stands in contrast with the preemptive problem, where customers can be moved instantaneously between the buffer and the stations that offer them service.

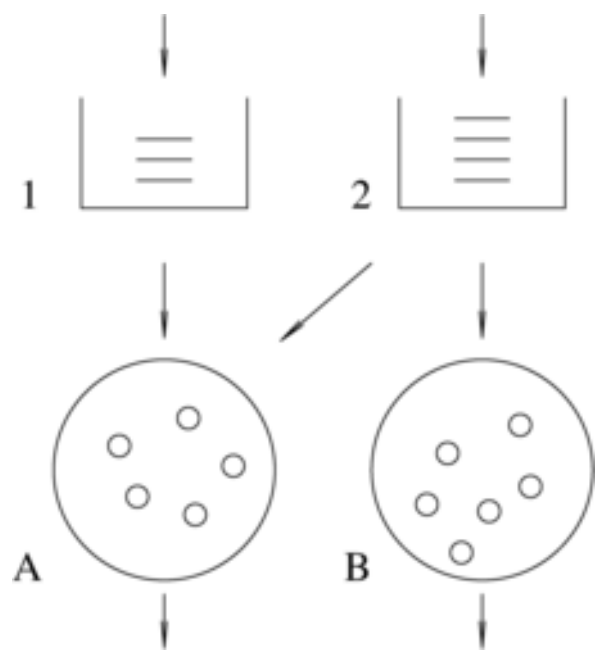

FIG. 2. *An example with three activities.*

We next argue that joint work conservation is in many cases a desired property for problems with preemption (for problems without preemption, work conservation is typically not optimal, cf. [3], but under appropriate conditions the distinction between the two disappears in the limit, as seen in our main result). We do not attempt a rigorous treatment, and instead describe an argument towards optimality in a particular case. We comment that even in the context of one station, proving that work conservation is optimal is not trivial (see [3]). However, it is believed that the argument below can be greatly generalized.

Example 2. The example consists of a Markovian system with two customer classes, two server types and no abandonment (see Figure 3). In (a) there is a customer of class 1 in the queue and a free server of type $B$. Although the free server cannot serve this customer, a rearrangement is possible so as to allow all customers to be served (b). We shall argue that "good" preemptive policies prefer option (b) over (a). Clearly this question must be coupled with the cost criteria. For concreteness, consider a cost of the form $\sum c_iY_i$ per unit time (weighted sum of queue lengths). Use the fact that optimality is obtained by feedback policies, that observe only the state: number of customers at each class present in the system. Given a feedback policy $\pi$ that leaves a customer in queue 1 when there is a free server in station $B$, we show there exists a policy $\hat\pi$ that pays a smaller cost on average. The argument is valid if $\mu_{2A} \le \max(\mu_{1A}, \mu_{2B})$. Arguing by coupling and assuming the system starts at the described state as time zero, let $\hat{\pi}$ move a class-2 customer, $C_2$, from station $A$ to station $B$ and let it move the class-1 customer waiting in the queue, $C_1$, into station $A$, until, at time $\tau$, there is an arrival or a service in one of the systems (i.e., the system under $\pi$ and that under $\hat{\pi}$). If the first to occur is an arrival or a service to a customer other than $C_1$ or $C_2$, $\hat{\pi}$ switches back to act like $\pi$ for all times. One can perform the coupling in such a way that service to $C_1$ under $\pi$ always occurs later than the first between service to $C_1$ and service to $C_2$ under $\hat{\pi}$. If either service to $C_1$ or to $C_2$ under $\hat{\pi}$ is the first to occur, $\hat{\pi}$ then mimics $\pi$ except for the single customer that under $\hat{\pi}$ is not present and may still be present for a while under $\pi$. In all cases the cost paid by $\hat{\pi}$ is not greater than that paid by $\pi$.

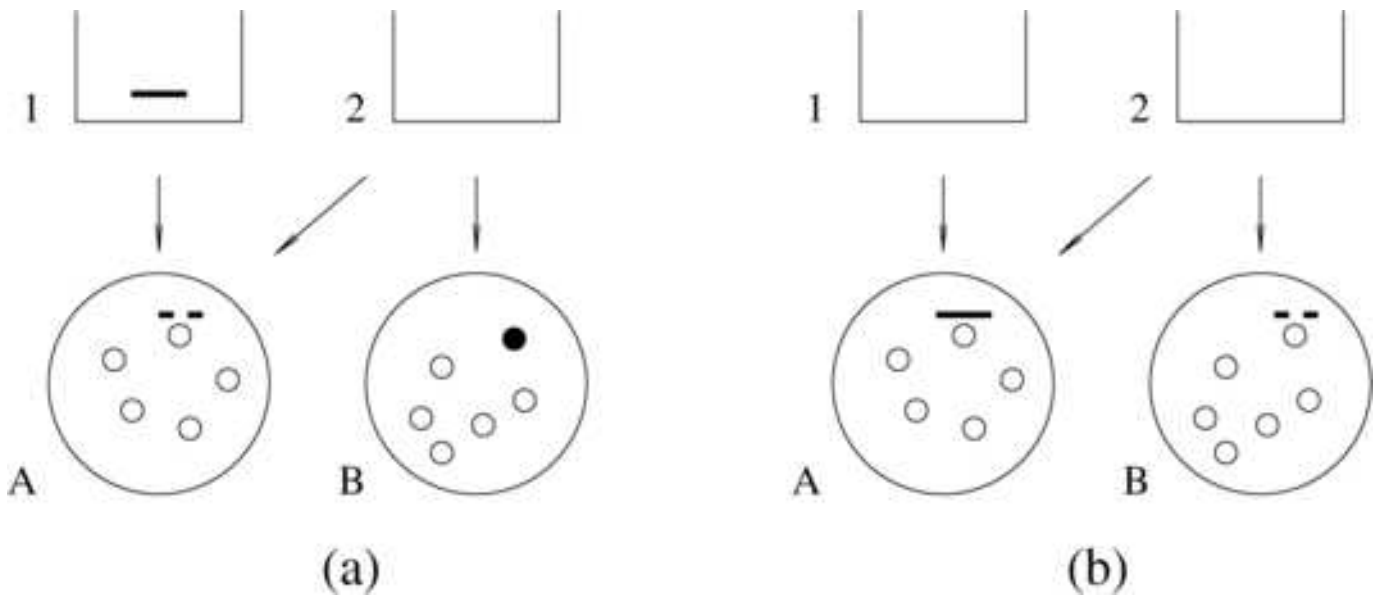

FIG. 3. (a) *A class-*1 *customer waiting in the queue, no free servers in station $A$ and one free server in station $B$.* (b) *A class-*2 *customer moved from station $A$ to station $B$ and the class-*1 *customer moved from the queue to station $A$.*

## 3. Proofs.

3.1. *Preliminary results.* In Section 2 we obtained the convenient representation (47) of the controlled process from (38)–(41). We begin by showing how an equation analogous to (47) is obtained for the prelimit model. To this end, let

$$M^n := e \cdot Y^n \wedge e \cdot Z^n \ge 0. \tag{61}$$

Denote also $\hat{M}^n = n^{-1/2} M^n$. By (27) and (28) we can write $e \cdot \hat{X}^n = e \cdot \hat{Y}^n - e \cdot \hat{Z}^n$ and therefore by (61),

$$e \cdot \hat{Y}^n = \hat{M}^n + (e \cdot \hat{X}^n)^+, \qquad e \cdot \hat{Z}^n = \hat{M}^n + (e \cdot \hat{X}^n)^-. \tag{62}$$

Let $u^n = \hat{Y}^n / (e \cdot \hat{Y}^n)$ when $e \cdot \hat{Y}^n > 0$, and set $u^n = e_1$ otherwise. Similarly, $v^n = \hat{Z}^n / (e \cdot \hat{Z}_n)$ if $e \cdot \hat{Z}^n > 0$ and otherwise set $v^n = e_{I+1}$. Letting $U^n =$ $(u^n, v^n)$, noting that $U^n$ takes values in $\mathbb{U}$, and using Lemma 1, (27), (28), (45) and linearity of $G$ on $D_G$, we have

$$\begin{aligned}\hat{Y}^n &= (\hat{M}^n + (e \cdot \hat{X}^n)^+)u^n,\\ \hat{Z}^n &= (\hat{M}^n + (e \cdot \hat{X}^n)^-)v^n\end{aligned} \tag{63}$$

and

$$\begin{aligned}\hat{\Psi}^n &= G(\hat{X}^n - \hat{Y}^n, -\hat{Z}^n)\\ &= G(\hat{X}^n - (e \cdot \hat{X}^n)^+u^n, -(e \cdot \hat{X}^n)^-v^n) - \hat{M}^n G(u^n, v^n)\\ &= \hat{G}(\hat{X}^n, U^n) - \hat{M}^n G(u^n, v^n).\end{aligned} \tag{64}$$

Defining

$$\begin{aligned}b_i^n(X, U) &= -\sum_j \mu_{ij}^n \hat{G}_{ij}(X, U) - \theta_i^n (e \cdot X)^+ u_i + \ell_i^n,\\ \check{G}_i^n(u, v) &= \sum_j \mu_{ij}^n G_{ij}(u, v) - \theta_i^n u_i,\end{aligned}$$

we obtain from (30)

$$\begin{aligned}\hat{X}_i^n(t) &= \hat{X}_i^{0,n} + r_i \hat{W}_i^n(t)\\ &\quad + \int_0^t \Bigg[\ell_i^n - \sum_j \mu_{ij}^n (\hat{G}_{ij}(\hat{X}^n(t), U_n(t)) - \hat{M}^n(t) G_{ij}(u^n(t), v^n(t)))\\ &\qquad\qquad - \theta_i^n (\hat{M}^n(t) u_i^n(t) + (e \cdot X^n(t))^+ u_i^n(t))\Bigg]\, ds,\end{aligned}$$

hence

$$\begin{aligned}\hat{X}^n(t) &= \hat{X}^{0,n} + r\hat{W}^n(t) + \int_0^t b^n(\hat{X}^n(s), U^n(s))\, ds\\ &\quad + \int_0^t \hat{M}^n(s) \check{G}^n(u^n(s), v^n(s))\, ds.\end{aligned} \tag{65}$$

It is useful to note that

$$\begin{aligned}&\|b^n(x, U) - b^n(y, U)\| \le c_1 \|x - y\|,\\ &\qquad \|b^n(x, U)\| \le c_1(1 + \|x\|), \qquad n \in \mathbb{N}, U \in \mathbb{U}, x, y \in \mathbb{R}^I,\end{aligned} \tag{66}$$

where $c_1$ does not depend on $n, U, x$ and $y$.

Part (i) of the following result proves (37). Part (ii) regards the proposed P-SCP.

LEMMA 3. *Fix $t$. Assume that $X^n(t) \in \mathbb{Z}_+^I$ satisfies the inequality $\|X^n(t) - nx^*\| \le \alpha_0 n$, where $\alpha_0 > 0$ is the constant from (36). Then the following conclusions hold.*

(i) *$X^n(t) \in \mathcal{X}^n$. Namely, there exist quantities $Y^n(t)$, $Z^n(t)$, $\Psi^n(t)$ satisfying (2)–(5) as well as (35).*

(ii) *Let $Y^n(t)$ and $Z^n(t)$ be given, define*

$$\Psi^n(t) = G(X^n(t) - Y^n(t), N^n - Z^n(t)) \tag{67}$$

*and assume that (4) and (35) are satisfied. Then (2), (3) and (5) also hold.*

REMARK. The lemma has two implications.

(a) Consider a jointly work-conserving P-SCP. If the inequality $\|X^n(t) - nx^*\| \le \alpha_0 n$ holds, then, by Lemma 3(i), $X^n(t) \in \mathcal{X}^n$. By definition of joint work conservation, a rearrangement of customers will be chosen so as to meet (35), and it follows that $M^n(t) = 0$. Now, let $\tau^n = \inf\{s : M^n(s) > 0\}$. The above discussion shows that if $\tau^n \le t$, then $\|X^n - nx^*\|_t^* \ge \alpha_0 n$. In particular, $\|X^n - nx^*\|_{\tau^n}^* \ge \alpha_0 n$ on $\{\tau^n < \infty\}$. This observation will be used in the next subsection.

(b) $p^*$ is a legitimate sequence of jointly work-conserving P-SCPs. As argued in Section 2, one only has to show that when $\|X^n(t) - nx^*\| \le \alpha_0 n$, $\Psi_{ij}^n(t) \ge 0$ holds. This is shown in Lemma 3(ii).

PROOF OF LEMMA 3. Consider part (ii) first. By the definition of $G$ (cf. Lemma 1), (2) and (3) hold. It remains to show that $\Psi_{ij}^n(t) \ge 0$ for all $(i,j) \in \mathcal{E}$. By Lemma 1 and (17),

$$\psi^* = G(x^*, \nu). \tag{68}$$

Since it is assumed that $M^n(t) = 0$ [cf. (61)], (62) implies that $\|Y^n(t)\| \vee \|Z^n(t)\| \le \|X^n(t) - nx^*\|$. Hence by linearity of the map $G$ on the domain $D_G$ and by (20), (44) and (36),

$$\begin{aligned}
\Psi_{ij}^n(t) &= G_{ij}(nx^*, n\nu) + G_{ij}(X^n(t) - nx^* - Y^n(t), N^n - n\nu - Z^n(t)) \\
&\ge n\psi_{ij}^* - C_G(\|X^n(t) - nx^* - Y^n(t)\| \vee \|N^n - n\nu - Z^n(t)\|) \\
&\ge n\psi_{ij}^* - 2C_G\|X^n(t) - nx^*\| - C_G\|N^n - n\nu\| \\
&\ge n\psi_{ij}^* - 2C_G\alpha_0 n - cn^{1/2} \\
&\ge 0,
\end{aligned}$$

where the last inequality holds for all $n$ large enough.

To prove part (i), simply set $Y^n(t) = (e \cdot X^n(t))^+ e_1$, $Z^n(t) = (e \cdot X^n(t))^- e_{J+1}$ and define $\Psi^n(t)$ via (67). The result now follows from part (ii). $\square$

The following lemma refers to instantaneous routing through nonblocked activities in the construction of the N-SCPs $p'$ of Section 2, by showing that one can find $\hat{\Psi}^n$ meeting (58).

LEMMA 4. *Let* $(\tilde{\Psi}, \tilde{X}, \tilde{Y}, \tilde{Z})$ *satisfy*

$$\begin{aligned} \sum_j \tilde{\Psi}_{ij} &= \tilde{X}_i - \tilde{Y}_i, && i \in \mathcal{I}, \\ \sum_i \tilde{\Psi}_{ij} &= -\tilde{Z}_j, && j \in \mathcal{J}, \\ \tilde{\Psi}_{ij} &= 0, && i \not\sim j, \end{aligned} \tag{69}$$

$$\tilde{Y}_i \geq 0, \qquad \tilde{Z}_j \geq 0, \qquad i \in \mathcal{I}, j \in \mathcal{J}. \tag{70}$$

*Assume that all components of* $\tilde{X}, \tilde{Y}, \tilde{Z}$ *and* $\tilde{\Psi}$ *are integers. Let a subset* $\mathcal{E}_1 \subset \mathcal{E}$ *("nonblocked" activities) be given. Then one can find* $(\Psi, X, Y, Z)$ *satisfying relations analogous to* (69), (70), *and*

$$X = \tilde{X}, \qquad Y \leq \tilde{Y}, \qquad Z \leq \tilde{Z}, \qquad \Psi \geq \tilde{\Psi}, \tag{71}$$

$$(i, j) \in \mathcal{E}_1 \quad \text{implies} \quad Y_i \wedge Z_j = 0. \tag{72}$$

PROOF. Define inductively a sequence $(X^{(k)}, Y^{(k)}, Z^{(k)}, \Psi^{(k)})$, $k = 0, \dots, k_1$, as follows. Let

$$(X^{(0)}, Y^{(0)}, Z^{(0)}, \Psi^{(0)}) = (\tilde{\Psi}, \tilde{X}, \tilde{Y}, \tilde{Z}).$$

Let $k \geq 0$ be given, for which $(X^{(k)}, Y^{(k)}, Z^{(k)}, \Psi^{(k)})$ is defined. If $Y^{(k)}$, $Z^{(k)}$ satisfy (72), then terminate, declaring $k_1 = k$. Otherwise, define $(X^{(k+1)}, Y^{(k+1)}, Z^{(k+1)}, \Psi^{(k+1)})$ as follows. Let $i_0$ be the smallest $i \in \mathcal{I}$ such that there is $j$ with $(i, j) \in \mathcal{E}_1$ and $Y_i^{(k)} \wedge Z_j^{(k)} > 0$. Let $j_0$ be the smallest such $j$. For $i \in \mathcal{I}, j \in \mathcal{J}$ define

$$\begin{aligned} X^{(k+1)} &= X^{(k)}, & Y_i^{(k+1)} &= Y_i^{(k)} - \mathbb{1}_{\{i=i_0\}}, \\ Z_j^{(k+1)} &= Z_j^{(k)} - \mathbb{1}_{\{j=j_0\}}, & \Psi_{ij}^{(k+1)} &= \Psi_{ij}^{(k)} + \mathbb{1}_{\{(i,j)=(i_0,j_0)\}}. \end{aligned}$$

Since by construction $0 \leq e \cdot Y^{(k)} = e \cdot \tilde{Y} - k$, the procedure must terminate. Defining $(X, Y, Z, \Psi) = (X^{(k_1)}, Y^{(k_1)}, Z^{(k_1)}, \Psi^{(k_1)})$ completes the proof. □

The following lemma will be useful in analyzing the N-SCPs $p'$. In order to show that $\hat{\Psi}^n$ tracks closely $\check{\Psi}^n$, estimates are required on $\|\hat{\Psi}^n(t) - \check{\Psi}^n(t)\|$. The lemma shows that it suffices to estimate only $(\hat{\Psi}^n(t) - \check{\Psi}^n(t))^+$.

LEMMA 5. *Let $(\psi, x, y, z)$ satisfy*

$$\sum_j \psi_{ij} = x_i - y_i, \qquad i \in \mathcal{I},$$

$$\sum_i \psi_{ij} = -z_j, \qquad j \in \mathcal{J},$$

$$\psi_{ij} = 0, \qquad i \not\sim j,$$

*and let $(\check{\psi}, \check{x}, \check{y}, \check{z})$ satisfy analogous relations. In addition, assume*

$$\text{if } i \sim j \text{ and } \psi_{ij} < \check{\psi}_{ij} \qquad \text{then } y_i \wedge z_j = 0, \tag{73}$$

*and $\check{y}_i \geq 0$, $i \in \mathcal{I}$, $\check{z}_j \geq 0$, $j \in \mathcal{J}$. Then*

$$\sum_{i,j} |\psi_{ij} - \check{\psi}_{ij}| \leq c \sum_{i,j} (\psi_{ij} - \check{\psi}_{ij})^+ + c\|x - \check{x}\|,$$

*where $c$ does not depend on $\psi, x, y, z, \check{\psi}, \check{x}, \check{y}$ or $\check{z}$.*

PROOF. Let $\varepsilon$ be an upper bound on $\psi_{ij} - \check{\psi}_{ij}$ for all $i, j$, and on $|x_i - \check{x}_i|$ for all $i$. Let $j_0$ be such that $z_{j_0} = 0$. Then

$$\sum_i \psi_{ij_0} = \sum_i \check{\psi}_{ij_0} + \check{z}_{j_0} \geq \sum_i \check{\psi}_{ij_0},$$

and since $\psi_{ij_0} \leq \check{\psi}_{ij_0} + \varepsilon$, $\psi_{ij_0} - \check{\psi}_{ij_0} \geq -c\varepsilon$ for every $i \sim j_0$. Let $i_0$ be such that $y_{i_0} = 0$. Then

$$\sum_j \psi_{i_0 j} = x_{i_0} \geq \check{x}_{i_0} - \varepsilon = \sum_j \check{\psi}_{i_0 j} + \check{y}_{i_0} - \varepsilon \geq \sum_j \check{\psi}_{i_0 j} - \varepsilon.$$

Since $\psi_{i_0 j} \leq \check{\psi}_{i_0 j} + \varepsilon$ for every $j$, we have $\psi_{i_0 j} - \check{\psi}_{i_0 j} \geq -c\varepsilon$ for every $j \sim i_0$. Thus we have shown that $|\psi_{ij} - \check{\psi}_{ij}| \leq c\varepsilon$ for every $(i, j)$, $i \sim j$, with either $y_i = 0$ or $z_j = 0$. In view of (73), we have shown that $|\psi_{ij} - \check{\psi}_{ij}| \leq c\varepsilon$ for every $(i, j)$, $i \sim j$, such that $\psi_{ij} < \check{\psi}_{ij}$. On the other hand, if $\psi_{ij} \geq \check{\psi}_{ij}$, then simply $|\psi_{ij} - \check{\psi}_{ij}| = \psi_{ij} - \check{\psi}_{ij} \leq \varepsilon$ by assumption. $\square$

Denote

$$J_t^n = \|\hat{Y}_t^n - \check{Y}_t^n\| + \|\hat{Z}_t^n - \check{Z}_t^n\|, \tag{74}$$

$$Q_1^n(t) = \int_0^t b^n(\hat{X}_s^n, U_s^n)\,ds, \qquad Q_2^n(t) = \int_0^t e^{-\gamma s} L(\hat{X}_s^n, U_s^n)\,ds. \tag{75}$$

Throughout, let $p$, $p^*$, $p'$ and $\zeta$ be as in Theorem 2, and let $f$ denote the unique $C^2_{\text{pol}}$ [$C^2_b$ under Assumption 3(iii)] solution to (49). Recall that $\mathscr{I} f = \int_0^\cdot f(s)\,ds$.

Proposition 1. *Items* (i)–(iii) *below hold under* $p$ *(in particular, under* $p^*$*) and under* $p'$.

(i) $(\bar X^n, \bar Y^n, \bar Z^n, \bar\Psi^n) \Rightarrow (x^*, 0, 0, \psi^*)$.

(ii) $(\hat W^n, \mathscr{I}\hat M^n, \hat X^n, Q_1^n, Q_2^n)$ *is tight.*

(iii) $(\hat W^n, \mathscr{I}\hat M^n) \Rightarrow (W, 0)$, *where* $W$ *is a standard Brownian motion on* $\mathbb{R}^I$. *Moreover, the following holds under* $p^*$ *and under* $p'$:

(iv) $|J^n|^*_{s,t} \to 0$ *and* $|\hat M^n|^*_{s,t} \to 0$ *in distribution, for every* $0 < s < t < \infty$.

For the proof see Section 3.2.

Lemma 6. *Under* $p$ *and under* $p'$ *one has the following. Denote by* $(X, Q_1, Q_2, W)$ *a limit point of* $(\hat X^n, Q_1^n, Q_2^n, \hat W^n)$ *along a subsequence. Let* $(F_t)$ *denote the filtration generated by* $(X, Q_1, W)$. *Then* $W$ *is an* $(F_t)$*-standard Brownian motion,* $X$, $Q_1$ *and* $Q_2$ *have continuous sample paths, and* $Q_1$ *has sample paths of bounded variation over finite time intervals. Moreover,* $\int e^{-\gamma s} Df(\hat X_s^n)\cdot dQ_1^n(s) \Rightarrow \int e^{-\gamma s} Df(X_s)\cdot dQ_1(s)$ *along the subsequence, where* $f$ *is the solution to* (49).

Proof. Based on Proposition 1, the proof of Lemma 6 is identical to that of Lemma 6 of [3] and is therefore omitted. □

Proposition 2. (i) *Assume either case* (i) *or* (ii) *of Assumption* 3 *holds. For either* $q = p$, *with fixed* $m_0 \in (m_L, m_A/2)$, *or for* $q = p'$, *with fixed* $m_0 \in (m_L, m_A(m_A - 2)(5m_A - 2)^{-1})$, *one has*

$$E_\zeta^q[(\|\hat X^n\|_t^*)^{m_0}] \le C(1+t)^{m_1} \tag{76}$$

*where* $C$, $m_1$ *do not depend on* $n$ *and* $t$. *In addition, the same conclusion holds in the case* $q = p'(\{h_n\})$, *for any sequence* $\{h_n\}$, *and the constants* $C$ *and* $m$ *do not depend on the sequence* $\{h_n\}$.

(ii) *Let case* (iii) *of Assumption* 3 *hold. Then*

$$E_\zeta^q[(\|\hat X^n\|_t^*)^{m_0}] \le \bar C e^{\bar C t}, \tag{77}$$

*for* $q = p'(\{h_n\})$, $m_0 \in (0, m_A(m_A - 2)(5m_A - 2)^{-1})$, *and* $\bar C$ *not depending on* $n$, $t$ *and the sequence* $\{h_n\}$.

In part (ii) of the proposition, the significance is not in the precise form of the upper bound but in the uniformity with respect to $\{h_n\}$.

For the proof of Proposition 2 see Section 3.3.

The method of [3], that we adopt here, is based on estimating the process

$$K_t^n = b(\hat X_t^n, U_t^n)\cdot Df(\hat X_t^n) + L(\hat X_t^n, U_t^n) - H(\hat X_t^n, Df(\hat X_t^n)) \ge 0, \tag{78}$$

where the inequality above follows from (50).

LEMMA 7. *Let the assumptions of Theorem 2 hold. For every sequence $p$ of admissible jointly work-conserving P-SCPs*

$$\liminf_{n\to\infty} E_\zeta^p \int_0^\infty e^{-\gamma t} L(\hat{X}_t^n, U_t^n)\, dt \geq V(x).$$

*Moreover, if for some $\bar{\delta} > 0$, $T > 0$ and $k \in (0, \infty]$, $q$ is an admissible SCP under which*

$$\lim_{n\to\infty} P_\zeta^q \left( \Omega_{n,k,T} \cap \left\{ \int_0^T e^{-\gamma t} K_t^n dt > \bar{\delta} \right\} \right) = 0, \tag{79}$$

*$\Omega_{n,k,T}$ being any event on which $\|\hat{X}^n\|_T^* \leq k$, then*

$$\limsup_{n\to\infty} E_\zeta^q \left[ \mathbb{1}_{\Omega_{n,k,T}} \int_0^T e^{-\gamma t} L(\hat{X}_t^n, U_t^n)\, dt \right] \leq V(x) + \bar{\delta}.$$

PROOF. Equipped with Proposition 1, Lemma 6 and Proposition 2 under Assumption 3(i) and (ii), and boundedness of $V$ in case of Assumption 3(iii), the proof is very similar to that of Theorem 4 of [3], and is therefore omitted. $\square$

PROOF OF THEOREM 2. Part (i): Established in Lemma 7.

Parts (ii) and (iii): We use Lemma 7 with $k = \infty$ and $\Omega_{n,k,T} = \Omega$. In view of this lemma it suffices to show that, for every $T > 0$, under both $p^*$ and $p'$, $\int_0^T e^{-\gamma s} K_s^n\, ds \to 0$ in probability. The function $h^*$ satisfies the following (see the proof of Theorem 1 of [2]):

$$H(x, Df(x)) = b(x, h^*(x)) \cdot Df(x) + L(x, h^*(x)), \qquad x \in \mathbb{R}^I. \tag{80}$$

Combining (78) and (80),

$$\begin{aligned} K_t^n &= (b(\hat{X}_t^n, U_t^n) - b(\hat{X}_t^n, h^*(\hat{X}_t^n))) \cdot Df(\hat{X}_t^n) \\ &\quad + L(\hat{X}_t^n, U_t^n) - L(\hat{X}_t^n, h^*(\hat{X}_t^n)). \end{aligned} \tag{81}$$

By definition of $b$ (46), and by (74),

$$\|b(\hat{X}_t^n, U_t^n) - b(\hat{X}_t^n, h^*(\hat{X}_t^n))\| \leq c J_t^n. \tag{82}$$

By (48) and (53),

$$L(\hat{X}_t^n, U_t^n) - L(\hat{X}_t^n, h^*(\hat{X}_t^n)) = \tilde{L}(\hat{X}_t^n, \Psi_1^n(t)) - \tilde{L}(\hat{X}_t^n, \check{\Psi}^n(t)) \tag{83}$$

where, using (63),

$$\begin{aligned} \Psi_1^n(t) &:= G(\hat{X}_t^n - (e \cdot \hat{X}_t^n)^+ u_t^n, -(e \cdot \hat{X}_t^n)^- v_t^n) \\ &= G(\hat{X}_t^n - \hat{Y}_t^n + \hat{M}_t^n u_t^n, -\hat{Z}_t^n + \hat{M}_t^n v_t^n) \\ &= \hat{\Psi}_t^n + \hat{M}_t^n G(u_t^n, v_t^n). \end{aligned} \tag{84}$$

Note that

$$\|\hat{\Psi}_t^n - \check{\Psi}_t^n\| = \|G(\hat{X}_t^n - \hat{Y}_t^n, -\hat{Z}_t^n) - G(\hat{X}_t^n - \check{Y}_t^n, -\check{Z}_t^n)\| \le cJ_t^n. \tag{85}$$

By Assumption 2, $\tilde{L}$ is uniformly continuous on compacts, hence there are functions $\alpha_k(\delta)$ with $\lim_{\delta\to 0}\alpha_k(\delta)=0$ such that

$$\begin{aligned}&|\tilde{L}(x,\psi)-\tilde{L}(x,\psi')| \le \alpha_k(\delta),\\ &\qquad \|x\|\vee\|\psi\|\vee\|\psi'\|\le k, \qquad \|\psi-\psi'\|\le\delta.\end{aligned} \tag{86}$$

Combining (81)–(86), $\|\hat{X}_t^n\|\vee\|\hat{\Psi}_t^n\|\le k$ implies

$$K_t^n \le cJ_t^n\|Df(\hat{X}_t^n)\| + \alpha_{ck}(J_t^n+\hat{M}_t^n) \le cJ_t^n\beta_k + \alpha_{ck}(J_t^n+\hat{M}_t^n), \tag{87}$$

where $\beta_k$ depends only on $k$. Since by (52), $e\cdot\check{Y}^n\wedge e\cdot\check{Z}^n=0$, we have

$$\hat{M}^n = e\cdot\hat{Y}^n\wedge e\cdot\hat{Z}^n \le |e\cdot\hat{Y}^n - e\cdot\check{Y}^n| + |e\cdot\hat{Z}^n - e\cdot\check{Z}^n| \le J^n. \tag{88}$$

Moreover, by (27), (28) and (52),

$$J^n \le c_1(\|\hat{X}^n\|+\|\hat{\Psi}^n\|), \tag{89}$$

for $c_1$ not depending on $n$. Fix $T$ and let $\Omega^{n,k,\varepsilon,\delta}$ denote the event that $\|\hat{X}^n\|_T^*\vee\|\hat{\Psi}^n\|_T^*\le k$ and $|J^n+\hat{M}^n|_{\varepsilon,T}^*\le\delta$. By Proposition 1,

$$\lim_k \lim_{\delta\to 0+}\liminf_{\varepsilon\to 0+}\liminf_n P^q(\Omega^{n,k,\varepsilon,\delta})=1, \tag{90}$$

for $q=p^*$ and for $q=p'$. Combining (87)–(89), on $\Omega^{n,k,\varepsilon,\delta}$ we have

$$\begin{aligned}0 &\le \int_0^T e^{-\gamma t}K_t^n\,dt\\ &\le c_2\varepsilon(k\beta_k+\alpha_{c_2k}(c_2k)) + c_2T(\beta_k\delta+\alpha_{c_2k}(\delta))\\ &=: \bar{\alpha}(k,\varepsilon,\delta),\end{aligned} \tag{91}$$

where $c_2$ does not depend on $n,\varepsilon,\delta,k$. Taking $n\to\infty$, then $\varepsilon\to 0$, $\delta\to 0$ and finally $k\to\infty$, using (90), it follows that $\int_0^T e^{-\gamma t}K_t^n\,dt\to 0$ in probability. Since $T$ is arbitrary, the result follows.

Part (iv): Recall that the sequence $p''=p'(\{h_n\})$ is defined like the sequence $p'$, except that for each $n$, $h$ is replaced by $h_n$. It suffices to show that for each $\bar{\delta}>0$ one can find $h^{\bar{\delta}}$ such that

$$\limsup_{n\to\infty} E_\zeta^q\int_0^\infty e^{-\gamma t}L(\hat{X}_t^n,U_t^n)\,dt \le V(x)+4\bar{\delta}, \tag{92}$$

where under $q$, $h_n$ are all equal to $h^{\bar{\delta}}$, because one can then find an appropriate sequence $\bar{\delta}_n\to 0$ and set $h_n=h^{\bar{\delta}_n}$ so as to obtain $\overline{V}(\zeta,p'')\le V(x)$. We will prove (92) by finding, for each $\bar{\delta}>0$, $h^{\bar{\delta}}$ and $T$ such that

$$\limsup_n E_\zeta^q\int_0^T e^{-\gamma t}L(\hat{X}_t^n,U_t^n)\,dt \le V(x)+3\bar{\delta} \tag{93}$$

and

$$\sup_n E^q_\zeta \int_T^\infty e^{-\gamma t} L(\hat{X}^n_t, U^n_t)\, dt \le \bar{\delta}. \tag{94}$$

Let $\bar{\delta} > 0$ be given. If $L$ is bounded, (94) holds for $T$ large. Otherwise, Assumption 3(i) or (ii) is in force, and (94) follows from Proposition 2(i) and Assumption 2(iv), for $T$ large enough. In particular, by the last statement of Proposition 2(i), such $T$ can be chosen independently of the sequence $\{h_n\}$. Let such $T$ be fixed. Using the two parts of Proposition 2(ii) for the different cases of Assumption 3, one can find $k \in (0, \infty)$, not depending on $\{h_n\}$, so large that

$$P^q_\zeta(\|\hat{X}^n\|^*_T \ge k) \le \bar{\delta}. \tag{95}$$

In addition, in view of (90), for $k$ large enough,

$$\limsup_{\delta\to 0^+} \liminf_{\varepsilon\to 0^+} \liminf_n E^q_\zeta \left[ \mathbb{1}_{(\Omega^{n,k,\varepsilon,\delta})^c} \int_0^T e^{-\gamma t} L(\hat{X}^n_t, U^n_t)\, dt \right] \le \bar{\delta}. \tag{96}$$

Let such $k$ be fixed.

We argue along the lines of the proof of Theorem 2 of [3]. From (80), the function $h^*$ satisfies $\varphi(x, h^*(x)) = \inf_{U\in\mathbb{U}} \varphi(x, U) =: \varphi^*(x)$, where

$$\varphi(x, U) = b(x, U) \cdot Df(x) + L(x, U).$$

For each $\bar{\varepsilon} > 0$ let $d(x, y, \bar{\varepsilon})$ denote the Euclidean distance from $y$ to the boundary $\partial B(x, \bar{\varepsilon} I^{1/2})$. Let $B_{x,\bar{\varepsilon}} = \bar{\varepsilon}\mathbb{Z}^I \cap B(x, \bar{\varepsilon} I^{1/2})$ and set

$$\tilde{d}(x, y, \bar{\varepsilon}) = \begin{cases} \dfrac{d(x, y, \bar{\varepsilon})}{\sum_{y'\in B_{x,\bar{\varepsilon}}} d(x, y', \bar{\varepsilon})}, & y \in B_{x,\bar{\varepsilon}}, \\ 0, & \text{otherwise.} \end{cases}$$

Let

$$h_{\bar{\varepsilon}}(x) = \sum_{y\in\bar{\varepsilon}\mathbb{Z}^I} \tilde{d}(x, y, \bar{\varepsilon}) h^*(y).$$

It is elementary to check that $x \mapsto h_{\bar{\varepsilon}}(x)$ is Lipschitz. Since $U \mapsto L(x, U)$ is convex by assumption, and $U \mapsto b(x, U)$ affine, the map $U \mapsto \varphi(x, U)$ is convex. By Jensen's inequality and uniform continuity of $(x, U) \to \varphi(x, U)$ and $x \mapsto \varphi^*(x)$ on $B(0, k)$, for sufficiently small $\bar{\varepsilon} > 0$,

$$\begin{aligned} b(x, h_{\bar{\varepsilon}}(x)) \cdot Df(x) + L(x, h_{\bar{\varepsilon}}(x)) &= \varphi(x, h_{\bar{\varepsilon}}(x)) \\ &= \varphi\left( x, \sum_y \tilde{d}(x, y, \bar{\varepsilon}) h^*(y) \right) \\ &\le \sum_y \tilde{d}(x, y, \bar{\varepsilon}) \varphi(x, h^*(y)) \end{aligned}$$

$$\le \sum_y \tilde{d}(x,y,\bar{\varepsilon})\varphi(y,h^*(y)) + \bar{\delta}/2$$
$$= \sum_y \tilde{d}(x,y,\bar{\varepsilon})\varphi^*(y) + \bar{\delta}/2$$
$$\le \varphi^*(x) + \bar{\delta} = H(x, Df(x)) + \bar{\delta},$$
$$x \in B(0,k).$$

Letting $h^{\bar{\delta}} = h_{\bar{\varepsilon}}$ we have [in place of (80)]

$$H(x, Df(x)) + \bar{\delta} \ge b(x, h^{\bar{\delta}}(x)) \cdot Df(x) + L(x, h^{\bar{\delta}}(x)), \qquad x \in B(0,k). \tag{97}$$

Hence, arguing as in part (iii) above, recalling that $h^{\bar{\delta}}$ is Lipschitz (a condition used in the proof of Proposition 1) one obtains in place of (91),

$$0 \le \int_0^T e^{-\gamma t} K_t^n \, dt \le \bar{\alpha}(k, \varepsilon, \delta) + \bar{\delta},$$

holding on $\Omega^{n,k,\varepsilon,\delta}$. Taking $\varepsilon$ and $\delta$ small enough so that $\bar{\alpha}(k,\varepsilon,\delta) \le \bar{\delta}$, $\int_0^T e^{-\gamma t} K_t^n \, dt \le 2\bar{\delta}$ on $\Omega^{n,k,\varepsilon,\delta}$. Thus $P(\Omega^{n,k,\varepsilon,\delta} \cap \{\int_0^T e^{-\gamma t} K_t^n \, dt > 2\bar{\delta}\})$ converges to zero as $n \to \infty$. By Lemma 7,

$$\limsup_{n\to\infty} E_\zeta^q \int_0^T e^{-\gamma t} L(\hat{X}_t^n, U_t^n) \, dt$$
$$\le V(x) + 2\bar{\delta} + \limsup_{n\to\infty} E_\zeta^q \biggl[ \mathbb{1}_{(\Omega^{n,k,\varepsilon,\delta})^c} \int_0^T e^{-\gamma t} L(\hat{X}_t^n, U_t^n) \, dt \biggr].$$

Using (96), we obtain (93). $\square$

LEMMA 8. *Let Assumption 4 hold. Then* $E(\|\hat{A}^n\|_t^*)^{m_A} \le c(1+t)^{m_A/2}$, *where* $c$ *does not depend on* $n$ *or* $t$.

For the proof see Lemma 2 of [3].

3.2. *Tightness estimates.* We prove Proposition 1. Most involved is the treatment of the nonpreemptive case. The main idea is a "bootstrap" argument (a variation of which is also used in the next subsection), where one first establishes tightness of the processes up to a certain stopping time, and then uses this to show that the probability that the stopping time is incurred in an arbitrary fixed time approaches zero. The proof is established in a number of steps.

Step 1. $n^{-1/2}\hat{W}^n \Rightarrow 0$.

Step 2. Under $p$, $\hat{M}^n \Rightarrow 0$ and $\bar{X}^n \Rightarrow x^*$. Under $p'$, $\bar{X}^n(\cdot \wedge \sigma^n) \Rightarrow x^*$, where $\sigma^n = \inf\{s > 0 : \mathscr{I}\hat{M}^n(s) \ge 1\} \wedge \Theta_n$.

Step 3. Under $p$, $(\bar{X}^n, \bar{Y}^n, \bar{Z}^n, \bar{\Psi}^n) \Rightarrow (x^*, 0, 0, \psi^*)$, and $(\hat{X}^n, \hat{W}^n, Q_1^n, Q_2^n)$ is tight.

Step 4. Under $p^*$, $J^n \Rightarrow 0$.

Step 5. Under $p'$, conclusions of step 3 hold, upon stopping all processes involved at $\sigma^n$. As a result, $\mathscr{I}\hat{M}^n \Rightarrow 0$, and conclusions analogous to those of step 3 hold under $p'$.

*Step* 1. Let $A_i$, $i \in \mathcal{I}$, $S_{ij}$, $i \sim j$, $R_i$, $i \in \mathcal{I}$, be independent Brownian motions with zero mean and variance given by $EA_i^2(1) = \lambda_i C_{U,i}^2$, $ES_{ij}^2(1) = \mu_{ij}$, $ER_i^2(1) = \theta_i$. Set $S_{ij} = 0$ for $i \not\sim j$, $A = (A_i)$, $S = (S_{ij})$ and $R = (R_i)$. For the fact $(\hat{A}^n, \hat{S}^n, \hat{R}^n) \Rightarrow (A, S, R)$ see [3], Lemma 4(i). Note that by (3) and (4), $\Psi_{ij}^n(t) \le N_j$ for every $i, j$ and $t$. Moreover, $\bar{Y}_i^n(t) \le e \cdot \bar{Y}_i^n(t) \le n^{-1}[e \cdot X^n(0) + \sum_i A_i^n(t)] \le c + n^{-1} \sum_i A_i^n(t) =: \zeta^n(t)$, where $c$ does not depend on $n$ and $t$. Hence by (31),

$$\|\hat{W}^n\|_t^* \le c_1 \|\hat{A}^n\|_t^* + c_1 \|\hat{S}^n\|_{c_1 t}^* + c_1 \|\hat{R}^n\|_{\zeta^n(t)}^*, \tag{98}$$

for a constant $c_1$. From this it is elementary to show that

$$n^{-1/2} \|\hat{W}^n\|_t^* \to 0 \qquad \text{in distribution, as } n \to \infty, t \ge 0. \tag{99}$$

*Step* 2. We show first that under $p$, $\hat{M}^n \Rightarrow 0$. Let $\tau^n = \inf\{s : \hat{M}_s^n > 0\}$. We shall show that, for every $T$, $P(\tau^n \le T) \to 0$ as $n \to \infty$; this implies that, for every $T$, $\|\hat{M}^n\|_T^* \to 0$ in distribution, and as a result $\hat{M}^n \Rightarrow 0$. Indeed, by (65),

$$\hat{X}^n(t \wedge \tau^n) = \hat{X}^{0,n} + \int_0^{t \wedge \tau^n} b^n(\hat{X}^n(s), U^n(s))\, ds + r\hat{W}^n(t \wedge \tau^n).$$

By (66),

$$\|\hat{X}^n(t \wedge \tau^n)\| \le c_2 + c_2 \int_0^{t \wedge \tau^n} (1 + \|\hat{X}^n(s)\|)\, ds + c_2 \|\hat{W}^n(t \wedge \tau^n)\|,$$

and it follows from Gronwall's inequality that

$$\|\bar{X}^n - x^*\|_{t \wedge \tau^n}^* = n^{-1/2} \|\hat{X}^n\|_{t \wedge \tau^n}^* \le c_3 n^{-1/2} e^{c_3 t} \|\hat{W}^n\|_t^*. \tag{100}$$

Hence by Remark (b) following Lemma 3 and (99),

$$P(\tau^n \le T) \le P(\|\bar{X}^n - x^*\|_{T \wedge \tau^n}^* > \alpha_0) \to 0. \tag{101}$$

As a result, $\hat{M}^n \Rightarrow 0$. By (99), (100) and (101) it follows that for every $t$, $\|\bar{X}^n - x^*\|_t^* \to 0$ in distribution as $n \to \infty$. As a result $\bar{X}^n \Rightarrow x^*$.

Next, under $p'$, recall that $\Theta_n$ is defined in (59) and let

$$\sigma^n = \inf\{t > 0 : \mathscr{I}\hat{M}^n(t) \ge 1\} \wedge \Theta_n. \tag{102}$$

By (65) and (66),

$$\|\hat{X}^n(t\wedge\sigma^n)\| \le c_4 + c_4\int_0^{t\wedge\sigma^n}(1+\|\hat{X}^n(s)\|)\,ds + c_4\|\hat{W}^n(t\wedge\sigma^n)\|.$$

Using again Gronwall's lemma and (99), we have

$$\|\bar{X}^n - nx^*\|^*_{t\wedge\sigma^n} \to 0 \qquad \text{in distribution, } t>0. \tag{103}$$

*Step* 3. This step refers to $p$ only. By (62),

$$e\cdot\bar{Y}^n = n^{-1/2}(\hat{M}^n + (e\cdot\hat{X}^n)^+) \le n^{-1/2}\hat{M}^n + \|\bar{X}^n - x^*\|.$$

Since $\bar{Y}_i \ge 0$, $i\in\mathcal{I}$, it follows that $\bar{Y}^n \Rightarrow 0$. By a similar argument, $\bar{Z}^n \Rightarrow 0$. By (67), (68) and linearity of the map $G$ on $D_G$,

$$\begin{aligned}\bar{\Psi}^n &= n^{-1}G(X^n - Y^n, N^n - Z^n)\\ &= G(\bar{X}^n - \bar{Y}^n, n^{-1}N^n - \bar{Z}^n)\\ &= G(\bar{X}^n - x^* - \bar{Y}^n, n^{-1}N^n - \nu - \bar{Z}^n) + G(x^*,\nu)\\ &= G(\bar{X}^n - x^* - \bar{Y}^n, n^{-1}N^n - \nu - \bar{Z}^n) + \psi^*.\end{aligned}$$

Since $\bar{Y}^n, \bar{Z}^n, (\bar{X}^n - x^*) \Rightarrow 0$ and by (20) and continuity of $G$, we obtain $\bar{\Psi}^n \Rightarrow \psi^*$.

We have now shown that $\|\bar{X}^n - x^*\|^*_t + \|\bar{Y}^n\|^*_t + \|\bar{Z}^n\|^*_t + \|\bar{\Psi}^n - \psi^*\|^*_t$ converges to zero in distribution, for every $t$. Hence $(\bar{X}^n, \bar{Y}^n, \bar{Z}^n, \bar{\Psi}^n) \Rightarrow (x^*, 0, 0, \psi^*)$.

Next we show that the sequence $(\hat{X}^n, \hat{W}^n, Q_1^n, Q_2^n)$ is tight in $(\mathbb{D}(\mathbb{R}^k))^3\times\mathbb{D}(\mathbb{R})$. We have shown already that $\hat{S}^n \Rightarrow S$, $\hat{R}^n \Rightarrow R$, $\bar{\Psi}^n \Rightarrow \psi^*$ and $\bar{Y}^n \Rightarrow 0$. An application of the time change lemma [5] shows that $\hat{S}^n_{ij}(\int_0^\cdot \bar{\Psi}^n_{ij}(s)\,ds) \Rightarrow S_{ij}(\psi^*\cdot)$ and $\hat{R}^n_i(\int_0^\cdot \bar{Y}^n_i(s)\,ds) \Rightarrow 0$. By (31) and (33) it follows that $\hat{W}^n \Rightarrow W$, a standard Brownian motion in $\mathbb{R}^I$.

By (65), (66),

$$\|\hat{X}^n(t)\| \le \|\hat{X}^{0,n}\| + c_5\|\hat{W}^n(t)\| + c_5\mathscr{I}\hat{M}^n(t) + c_5\int_0^t(1+\|\hat{X}^n(s)\|)\,ds$$

and therefore by Gronwall's inequality,

$$\|\hat{X}^n\|^*_t \le c_6 e^{c_6 t}(1+\mathscr{I}\hat{M}^n_t + \|\hat{W}^n\|^*_t). \tag{104}$$

Using tightness of $\hat{W}^n$ and $\hat{M}^n$, it follows that for every $t$,

$$\lim_{m\to\infty}\limsup_{n\to\infty} P(\|\hat{X}^n\|^*_t \ge m) = 0. \tag{105}$$

Fix $T$. By (65) and (66), for $s,t\in[0,T]$, $s<t$,

$$\begin{aligned}\|\hat{X}^n(t) - \hat{X}^n(s)\| &\le c_7\|\hat{W}^n(t) - \hat{W}^n(s)\|\\ &\qquad + c_7(t-s)(1+|\hat{M}^n|^*_T + \|\hat{X}^n\|^*_T).\end{aligned} \tag{106}$$

Let $w(x,S)=\sup_{s,t\in S}\|x(s)-x(t)\|$ where $S\subset[0,T)$, and let

$$w'_T(x,\delta)=\inf\max_{1\le i\le v} w(x,[t_{i-1},t_i)),$$

where the infimum is over all decompositions $[t_{i-1},t_i)$, $1\le i\le v$, of $[0,T)$ such that $t_i-t_{i-1}>\delta$, $1\le i\le v$ (cf. [5], page 171). The notation

$$w_T(x,\delta)=\sup_{0\le s<t\le(s+\delta)\wedge T} w(x,[s,t]) \tag{107}$$

will also be useful. It follows from Theorem 16.8 of [5] and tightness of $\hat W^n$ that for each $t$ and $\varepsilon$, $\lim_{\delta\to0}\limsup_{n\to\infty}P(w'_t(\hat W^n,\delta)\ge\varepsilon)=0$. Hence by tightness of $\hat M^n$ and using (106), for each $t\le T$ and $\varepsilon$,

$$\lim_{\delta\to0}\limsup_{n\to\infty}P(w'_t(\hat X^n,\delta)\ge\varepsilon)=0. \tag{108}$$

Since $T$ is arbitrary, (105) and (108) imply tightness of $\hat X^n$, using Theorem 16.8 of [5]. By (66), (75) and Assumption 2, there is a constant $c_8$ such that

$$\|Q_1^n\|_t^*\vee|Q_2^n|_t^*\le c_8t(1+\|\hat X^n\|_t^*)^{m_L},$$

and for $s,t\in[0,T]$, $s<t$,

$$\|Q_1^n(t)-Q_1^n(s)\|\vee|Q_2^n(t)-Q_2^n(s)|\le c_8(t-s)(1+\|\hat X^n\|_T^*)^{m_L}.$$

Hence, using Theorem 16.8 of [5], tightness of $Q_1^n$ and $Q_2^n$ follows from (105). We have shown that $\hat X^n$, $\hat W^n$, $Q_1^n$ and $Q_2^n$ are tight, and using Proposition 3.2.4 of [6], it follows that $(\hat X^n,\hat W^n,Q_1^n,Q_2^n)$ is tight.

*Step* 4. Fix $T$. By (1) and (55), under $p^*$ one has $|J^n|_T^*\le2(I+J)n^{-1/2}$ on the event $\{\|\bar X^n-x^*\|_T^*\le\alpha_0\}$. Since $\|\bar X^n-x^*\|_T^*$ converges to zero in distribution, so does $|J^n|_T^*$, and since $T$ is arbitrary, $J^n\Rightarrow0$.

*Step* 5. This step refers to $p'$. Let $\sigma^n$ be as in (102) and recall (103). Reviewing step 3 shows that all its conclusions still hold under $p'$ in place of $p$, upon replacing $\bar X^n$ by $\bar X^n(\cdot\wedge\sigma^n)$, $\bar Y^n$ by $\bar Y^n(\cdot\wedge\sigma^n)$, and similar substitutions for the processes $\bar\Psi^n$, $\hat W^n$, $\hat X^n$, $Q_1^n$ and $Q_2^n$. As a result,

$$(\bar X^n_{\cdot\wedge\sigma^n},\bar Y^n_{\cdot\wedge\sigma^n},\bar Z^n_{\cdot\wedge\sigma^n},\bar\Psi^n_{\cdot\wedge\sigma^n})\Rightarrow(x^*,0,0,\psi^*), \tag{109}$$

$$(\hat X^n_{\cdot\wedge\sigma^n},\hat W^n_{\cdot\wedge\sigma^n},Q_1^n(\cdot\wedge\sigma^n),Q_2^n(\cdot\wedge\sigma^n))\qquad\text{is tight.} \tag{110}$$

Let $T$ be fixed and denote $T^n=T\wedge\sigma^n$. For $i\sim j$ let

$$\Lambda_{ij}^n=\hat\Psi_{ij}^n-\check\Psi_{ij}^n. \tag{111}$$

Let $\Omega^{n,k}$ denote the event

$$\{\|\hat{X}^n\|^*_{T^n} \vee \|\hat{Y}^n\|^*_{T^n} \vee \|\hat{Z}^n\|^*_{T^n} \vee \|\hat{\Psi}^n\|^*_{T^n} \le k\}.$$

By tightness of $\hat{X}^n(\cdot \wedge \sigma^n)$ and by (62), (64),

$$\lim_{k\to\infty} \liminf_{n\to\infty} P(\Omega^{n,k}) = 1. \tag{112}$$

We will show that

$$\lim_{\varepsilon_0\to 0+} \limsup_n P(\sigma^n \le \varepsilon_0) = 0, \tag{113}$$

and that for every $i \sim j$ and small enough $\varepsilon_0 > 0$,

$$\limsup_n P\biggl(\sigma^n > \varepsilon_0, \sup_{[\varepsilon_0,T^n]} \Lambda^n_{ij} > c\varepsilon\biggr) = 0. \tag{114}$$

To this end let $\varepsilon_0 \in (0,T)$. Fix $i,j$, $i \sim j$. If $\Lambda^n_{ij} > 0$ on $[s,r)$, no customers are routed on activity $(i,j)$ on this time interval and by (57) and (9),

$$\Psi^n_{ij}(t) = \Psi^n_{ij}(s) - \Delta^n_{ij}(s,t), \qquad t \in [s,r),$$

where

$$\Delta^n_{ij}(s,t) = S^n_{ij}(\mathscr{I}\Psi^n_{ij}(t)) - S^n_{ij}(\mathscr{I}\Psi^n_{ij}(s)). \tag{115}$$

Therefore

$$\hat{\Psi}^n_{ij}(t) = \hat{\Psi}^n_{ij}(s) - n^{-1/2}\Delta^n_{ij}(s,t). \tag{116}$$

Now,

$$P\biggl(\sigma^n > \varepsilon_0, \sup_{t\in[\varepsilon_0,T^n]} \Lambda^n_{ij}(t) > 3\varepsilon\biggr) \le P((\Omega^{n,k})^c) + P(\Omega^{n,k}_1) + P(\Omega^{n,k}_2), \tag{117}$$

where

$$\begin{aligned}\Omega^{n,k}_1 &= \Omega^{n,k} \cap \{\sigma^n > \varepsilon_0\} \\ &\qquad \cap \Bigl\{\exists\, 0 \le s \le r \le T^n : \Lambda^n_{ij}(s) \le \varepsilon, \inf_{t\in(s,r)} \Lambda^n_{ij}(t) \ge \varepsilon, \Lambda^n_{ij}(r) \ge 3\varepsilon\Bigr\}, \\ \Omega^{n,k}_2 &= \Omega^{n,k} \cap \{\sigma^n > \varepsilon_0\} \cap \Bigl\{\inf_{t\in[0,\varepsilon_0]} \Lambda^n_{ij} \ge \varepsilon\Bigr\}.\end{aligned}$$

Let $k$ be fixed. By the Hölder assumption on $h_1, h_2$ away from $e \cdot \hat{X}^n = 0$, there are positive constants $c'_k$, $\bar{c}_k$, $p_k$ (depending on $k$ and $\varepsilon$ but not on $n$) such that the following holds on $\Omega^{n,k}$:

$$\begin{aligned}|\check{\Psi}^n_{ij}(t) - \check{\Psi}^n_{ij}(s)| &= |G_{ij}(\hat{X}^n(t) - \check{Y}^n(t), -\check{Z}^n(t)) \\ &\qquad - G_{ij}(\hat{X}^n(s) - \check{Y}^n(s), -\check{Z}^n(s))|\end{aligned}$$

$$
\begin{aligned}
&\le c(\|(\hat X^n(t)-\check Y^n(t))-(\hat X^n(s)-\check Y^n(s))\|+\|\check Z^n(t)-\check Z^n(s)\|)\\
&\le c_k'(\|\hat X^n(t)-\hat X^n(s)\|+\|\hat X^n(t)-\hat X^n(s)\|^{p_k})\\
&\quad+\frac{\varepsilon}{4}\mathbb{1}_{\{|e\cdot\hat X^n_s|<\varepsilon/8\}}+\frac{\varepsilon}{4}\mathbb{1}_{\{|e\cdot\hat X^n_t|<\varepsilon/8\}}\\
&\le \bar c_k\|\hat X^n(t)-\hat X^n(s)\|^{p_k}+\frac{\varepsilon}{2}.
\end{aligned}\tag{118}
$$

By (23), for $0\le s\le t\le r\le T^n$ as above,

$$
\begin{aligned}
n^{-1/2}\Delta^n_{ij}(s,t)&=\hat S^n_{ij}(\mathscr{I}\bar\Psi^n_{ij}(t))-\hat S^n_{ij}(\mathscr{I}\bar\Psi^n_{ij}(s))\\
&\quad+n^{1/2}\mu^n_{ij}\int_s^t(\bar\Psi^n_{ij}(\vartheta)-\psi^*_{ij})\,d\vartheta\\
&\quad+n^{1/2}(\mu^n_{ij}-\mu_{ij})\psi^*_{ij}(t-s)+n^{1/2}\mu_{ij}\psi^*_{ij}(t-s).
\end{aligned}\tag{119}
$$

Since $\Psi^n_{ij}\le N^n_j$, it follows that $\bar\Psi^n_{ij}$ are all bounded uniformly by a constant. Since the jumps of the process $\hat S^n_{ij}$ are all of size $n^{-1/2}$, a use of (12.9) of [5] shows that $|\hat S^n_{ij}(\mathscr{I}\bar\Psi^n_{ij}(t))-\hat S^n_{ij}(\mathscr{I}\bar\Psi^n_{ij}(s))|\le 2w'_{c_9T^n}(\hat S^n_{ij},c_9(t-s))+n^{-1/2}$, where $c_9>0$ is a constant. Since $\|\bar\Psi^n-\psi^*\|^*_{T^n}\le kn^{-1/2}$ on $\Omega^{n,k}$, using (19), and assuming $n$ is large enough we therefore have

$$
n^{-1/2}\Delta^n_{ij}(s,t)\ge -2w'_{c_9T^n}(\hat S^n_{ij},c_9(t-s))-n^{-1/2}+c_{10}n^{1/2}(t-s),\tag{120}
$$

where $c_{10}=\mu_{ij}\psi^*_{ij}/2$. Hence on $\Omega^{n,k}_1$, noting that $\Lambda^n_{ij}(s)\le\varepsilon$,

$$
\begin{aligned}
\Lambda^n_{ij}(t)&=\Lambda^n_{ij}(s)+(\hat\Psi^n_{ij}(t)-\hat\Psi^n_{ij}(s))-(\check\Psi^n_{ij}(t)-\check\Psi^n_{ij}(s))\\
&\le -n^{-1/2}\Delta^n_{ij}(s,t)+\bar c_k\|\hat X^n(t)-\hat X^n(s)\|^{p_k}+2\varepsilon\\
&\le 2w'_{c_9T^n}(\hat S^n_{ij},c_9(t-s))-c_{10}n^{1/2}(t-s)+\bar c_k\|\hat X^n(t)-\hat X^n(s)\|^{p_k}+2\varepsilon.
\end{aligned}
$$

Since on $\Omega^{n,k}_1$, $\Lambda^n_{ij}(r)\ge 3\varepsilon$,

$$
\begin{aligned}
\varepsilon&\le 2w'_{c_9T^n}(\hat S^n_{ij},c_9(r-s))\\
&\quad-c_{10}n^{1/2}(r-s)+\bar c_k\|\hat X^n(r)-\hat X^n(s)\|^{p_k},
\end{aligned}\tag{121}
$$

$$
P(\Omega^{n,k}_1)\le P(\Omega^{n,k}_3)+P(\Omega^{n,k}_4),
$$

where

$$
\Omega^{n,k}_3=\{\exists\,0\le s\le r\le T^n: r-s\le an^{-1/2},\ (121)\text{ holds}\},
$$

$$
\Omega^{n,k}_4=\{\exists\,0\le s\le r\le T^n: r-s> an^{-1/2},\ (121)\text{ holds}\}.
$$

On $\Omega^{n,k}_3$,

$$
\varepsilon\le 2w'_{c_9T^n}(\hat S^n_{ij},c_9an^{-1/2})+\bar c_k w'_{T^n}(\hat X^n,an^{-1/2})^{p_k}
$$

and since $(\hat{S}^n, \hat{X}^n(\cdot \wedge \sigma^n))$ are tight, for every $a$,

$$\lim_n P(\Omega_3^{n,k}) = 0.$$

On $\Omega_4^{n,k}$,

$$c_{10} a \le 2\|\hat{S}_{ij}^n\|_{c_9 T^n}^* + \bar{c}_k (2\|\hat{X}^n\|_{T^n}^*)^{p_k}$$

and by tightness of $(\hat{S}^n, \hat{X}^n(\cdot \wedge \sigma^n))$,

$$\lim_{a\to\infty} \limsup_{n\to\infty} P(\Omega_4^{n,k}) = 0.$$

As a result we have, for every $k$,

$$\lim_{n\to\infty} P(\Omega_1^{n,k}) = 0. \tag{122}$$

Regarding $\Omega_2^{n,k}$, substituting $s = 0$ and $t \in [0, \varepsilon_0]$ in (118)–(120),

$$\Lambda_{ij}^n(t) \le \Lambda_{ij}^n(0) + 2\|\hat{S}_{ij}^n\|_{c_9\varepsilon_0}^* - c_{10} n^{1/2} t + \bar{c}_k \|\hat{X}^n(t) - \hat{X}^{0,n}\|^{p_k} + \varepsilon/2.$$

Since $\Lambda_{ij}^n(\varepsilon_0) \ge \varepsilon$, we have on $\Omega_2^{n,k}$

$$\varepsilon/2 \le \Lambda_{ij}^n(0) + 2\|\hat{S}_{ij}^n\|_{c_9\varepsilon_0}^* - c_{10} n^{1/2} \varepsilon_0 + \bar{c}_k \|\hat{X}^n(\varepsilon_0) - \hat{X}^{0,n}\|^{p_k}.$$

By tightness of the random variables $\xi^n := 2\|\hat{S}_{ij}^n\|_{c_9\varepsilon_0}^* + \bar{c}_k\|\hat{X}^n(\varepsilon_0 \wedge \sigma^n) - \hat{X}^{0,n}\|^{p_k}$, the fact that $\sigma^n > \varepsilon_0$ on $\Omega_2^{n,k}$, and since $\Lambda_{ij}^n(0)$ are bounded [as follows from (21), (22)],

$$\limsup_{n\to\infty} P(\Omega_2^{n,k}) \le \limsup_{n\to\infty} P(\xi^n \ge c\varepsilon_0 n^{1/2}) = 0. \tag{123}$$

Combining (117), (122), (123) and (112), we obtain (114).

Now, on $\Omega^{n,k}$, $\hat{M}^n(\sigma^n \wedge \varepsilon_0) \le e \cdot \hat{Y}^n(\sigma^n \wedge \varepsilon_0) \le k$. Hence by (102),

$$P(\sigma^n \le \varepsilon_0) \le P(\mathscr{I}\hat{M}^n(\sigma^n \wedge \varepsilon_0) \ge 1) \le P((\Omega^{n,k})^c) + \mathbb{1}_{\{\varepsilon_0 k \ge 1\}}.$$

Therefore, for every $k$,

$$\limsup_{\varepsilon_0\to 0+} \limsup_n P(\sigma^n \le \varepsilon_0) \le \limsup_n P((\Omega^{n,k})^c),$$

and (113) follows using (112).

Having established (113) and (114), we argue as follows. Lemma 5 and the property (58) of the policy $p'$ imply that

$$\|\Lambda^n(t)\| \le c_{11} \sum_{i\sim j} (\Lambda_{ij}^n(t))^+. \tag{124}$$

By (53) and the uniqueness statement in Lemma 1, $\hat{X}^n_i - \check{Y}^n_i = \sum_j \check{\Psi}^n_{ij}$ for every $i$. This and (27), along with an analogous argument for $\check{Z}^n$, imply that

$$\|\hat{Y}^n - \check{Y}^n\| + \|\hat{Z}^n - \check{Z}^n\| \le c_{12}\|\Lambda^n\|. \tag{125}$$

Combining (88), (124) and (125), on $\Omega^{n,k}$, for every $t \in (\varepsilon_0, T]$,

$$\mathscr{I}\hat{M}^n(t) \le c_{13}k\varepsilon_0 + c_{13}t \sum_{i\sim j} \sup_{[\varepsilon_0,t]} (\Lambda^n_{ij})^+. \tag{126}$$

Hence, with $b_0$ as in (60),

$$\begin{aligned} P(\sigma^n \le T) &\le P((\Omega^{n,k})^c) + P(\sigma^n \le \varepsilon_0) \\ &\qquad + P(\sigma^n > \varepsilon_0, \mathscr{I}\hat{M}^n(T^n) \ge 1) \\ &\qquad + P\Big(\sigma^n > \varepsilon_0, \max_{i\sim j} \sup_{[\varepsilon_0,T^n]} (\Lambda^n_{ij})^+ \ge b_0\Big) \\ &\le P((\Omega^{n,k})^c) + P(\sigma^n \le \varepsilon_0) \\ &\qquad + P\Big(\sigma^n > \varepsilon_0, c_{13}T\sum_{i\sim j} \sup_{[\varepsilon_0,T^n]} (\Lambda^n_{ij})^+ \ge 1 - c_{13}k\varepsilon_0\Big) \\ &\qquad + P\Big(\max_{i\sim j} \sup_{[\varepsilon_0,T^n]} (\Lambda^n_{ij})^+ \ge b_0\Big). \end{aligned} \tag{127}$$

Taking $\varepsilon_0$ small enough and using (112), (113) and (114), we have that

$$\lim_n P(\sigma^n \le T) = 0. \tag{128}$$

Since $T$ is arbitrary, we finally have from (109), (110) and (128) that

$$(\bar{X}^n, \bar{Y}^n, \bar{Z}^n, \bar{\Psi}^n) \Rightarrow (x^*, 0, 0, \psi^*), \qquad (\hat{X}^n, \hat{W}^n, Q^n_1, Q^n_2) \text{ is tight.}$$

Also, with (128), the relations (114) and (124) show that

$$\lim_n P\Big(\sup_{[\varepsilon_0,T]} \|\Lambda^n\| > c\varepsilon\Big) = 0. \tag{129}$$

In view of (74) and (125), we have from (129) that $|J^n|^*_{s,t}$ converges to zero in distribution, for every $0 < s < t < \infty$. As follows from (88), a similar statement holds for $|\hat{M}^n|^*_{s,t}$. Moreover, using again (126), now equipped with (129), letting $n \to \infty$, then $\varepsilon_0 \to 0^+$ and finally $k \to \infty$, we obtain that $\mathscr{I}\hat{M}^n(T) \to 0$ in distribution, and since $T$ is arbitrary, $\mathscr{I}\hat{M}^n \Rightarrow 0$. This completes the proof of Proposition 1.

3.3. *Large time estimates.* In this section we prove Proposition 2. The following estimate from [2] is used in a crucial way.

PROPOSITION 3. *Let Assumption* 3(i) *or* (ii) *hold. Let (27)–(30) hold. Then*

$$\|\hat{X}^n(t)\| \le C(1+t)^m(\|x\| + \|\hat{W}^n\|_t^* + |\hat{M}^n|_t^*), \tag{130}$$

*for $C$ and $m$ depending only on the model parameters, and in particular, not depending on $x$, $\hat{W}^n$ or on the SCP.*

PROOF. We use the following result from [2]: If

$$\begin{aligned} &x_i = w_i - \sum_j \mu_{ij}\mathscr{I}\psi_{ij} - \theta_i\mathscr{I}y_i, \\ &\psi = G(x-y,-z), \qquad e\cdot y \wedge e\cdot z = 0, \end{aligned} \tag{131}$$

then

$$\|x\| \le c(1+t)^m\|w\|_t^*, \tag{132}$$

where the constants $c,m$ do not depend on $\psi,x,y,z$: Under Assumption 3(i), (132) follows from the proof of Corollary 1 of [2]. Under Assumption 3(ii), (132) follows from Theorem 3 and Lemma 3 of [2].

Let

$$W_1^n(t) = \hat{W}^n(t) + r^{-1}\int_0^t \hat{M}^n(s)\check{G}^n(u^n(s),v^n(s))\,ds. \tag{133}$$

We have from (65) $\hat{X}_t^n = \hat{X}^{0,n} + \int_0^t b^n(\hat{X}_s^t, U_s^n)\,ds + rW_1^n(t)$. Lemma 2 implies that there are $\Psi_1^n$, $Y_1^n$ and $Z_1^n$ such that (38)–(41) hold. Hence $\hat{X}^n, Y_1^n, Z_1^n$, $\Psi_1^n, W_1^n$ satisfy relations analogous to (131). As a result, a relation as in (132) holds, and using (133) we obtain (130). □

REMARK. The only two places in this paper where Assumption 3(i), (ii) is used are in Theorem 1 and in obtaining (132) from (131). In fact, also the proof of Theorem 1 (carried out in [2]) uses this assumption only in order to establish (132), and if the implication "(131) implies (132)" holds true in greater generality, then so do Theorems 1 and 2.

PROOF OF PROPOSITION 2. Most of the proof deals with part (i). Let Assumption 3(i) or (ii) hold. First consider policies of the form $p$. By (98),

$$\|\hat{W}^n\|_t^* \le c_1(\|\hat{A}^n\|_t^* + \|\hat{S}^n\|_{c_1 t}^* + \|\hat{R}^n\|_{\zeta^n(t)}^*).$$

By Assumption 4, applying Lemma 8 shows that, for an appropriate constant $c_2$, $E(\|\hat{A}^n\|_t^*)^{m_A} \le c_2(1+t)^{c_2}$, and a similar estimate holds for $\hat{S}^n$.

Conditioning first on $\zeta^n$ and using the independence of $A^n$ and $R^n$, a similar estimate follows for the last term above. As a result,

$$E(\|\hat{W}^n\|_t^*)^{m_A} \le c_2(1+t)^{c_2}. \tag{134}$$

As in the proof of Proposition 1, let $\tau^n = \inf\{t : \hat{M}_t^n > 0\}$. Under $\{\tau^n > t\}$, $|\hat{M}^n|_t^* = 0$ and therefore by Proposition 3,

$$\|\hat{X}_t^n\|^* \le C(1+t)^m(\|x\| + \|\hat{W}^n\|_t^*), \qquad \{\tau^n > t\}. \tag{135}$$

On $\{\tau^n \le t\}$, by (3), (5) and (20),

$$|\hat{M}^n|_t^* = n^{-1/2}|M^n|_t^* \le n^{-1/2}|e \cdot Z^n|_t^* \le n^{-1/2} e \cdot N^n \le c_3 n^{1/2},$$

and by Proposition 3,

$$\|\hat{X}^n\|_t^* \le C(1+t)^m(\|x\| + \|\hat{W}^n\|_t^* + n^{1/2}), \qquad \{\tau^n \le t\}. \tag{136}$$

Combining (134), (135) and Lemma 3 [in particular, Remark (a) that follows],

$$\begin{aligned} P(\tau^n \le t) &\le P(\|\hat{X}^n\|_{t\wedge\tau^n}^* \ge \alpha_0 n^{1/2}) \\ &\le c n^{-m_A/2} E(\|\hat{X}^n\|_{t\wedge\tau^n}^*)^{m_A} \\ &\le c_4 n^{-m_A/2}(1+t)^{c_4}. \end{aligned} \tag{137}$$

Therefore by (136), (137) and the Hölder inequality, with $q^{-1} + \bar{q}^{-1} = 1$ and $q m_0 = m_A$,

$$\begin{aligned} E(\|\hat{X}^n\|_t^*)^{m_0} &\le E[(\|\hat{X}^n\|_t^*)^{m_0} \mathbb{1}_{\{\tau^n > t\}}] \\ &\quad + (E(\|\hat{X}^n\|_t^*)^{m_A})^{1/q}(P(\tau^n \le t))^{1/\bar{q}} \\ &\le c_5(1+t)^{c_5} + c_5 n^{(2q)^{-1}m_A} n^{-(2\bar{q})^{-1}m_A}(1+t)^{c_4/\bar{q}} \\ &\le c_6(1+t)^{c_6}, \end{aligned} \tag{138}$$

where the inequality $q^{-1} - \bar{q}^{-1} \le 0$, used on the last line above, follows from $m_0 < m_A/2$.

Next consider the policy $p'$. Let $b_1 := (2b_0) \vee 13$, where $b_0$ is the constant from (60). Let

$$\vartheta_n = \inf\Big\{t \ge 0 : \max_{i\sim j} \Lambda_{ij}^n \ge b_1\Big\}.$$

By (88), (124) and (125),

$$\hat{M}^n \le J^n \le c_7\|\Lambda^n\| \le c_8 \max_{i\sim j}(\Lambda_{ij}^n)^+.$$

Letting $T_n = T \wedge \vartheta_n$, it follows that

$$|\hat{M}^n|_{T_n}^* \le c_8 b_1.$$

Hence by Proposition 3, we have

$$\|\hat{X}^n\|^*_{T_n} \le C(1+T_n)^m(\|x\| + \|\hat{W}^n\|^*_{T_n} + c_8 b_1) \tag{139}$$
$$\le c_9(1+T)^m(1+\|\hat{W}^n\|^*_T),$$

where $c_9$ depends on $x$ but not on $T$ or $n$. We establish below the estimate

$$P(\vartheta_n \le T) \le c_0 n^{1/4-m_A/8}(1+T)^m, \tag{140}$$

where $c_0$ and $m$ are constants that do not depend on $T$ or $n$. Repeating the argument of (138), with $\vartheta_n$ in place of $\tau^n$, and (139) [resp., (140)] in place of (135) [resp., (137)], shows that $E(\|\hat{X}^n\|^*_T)^{m_0} \le c_{10}(1+T)^{c_{10}}$, with $qm_0 = m_A$, provided that $m_0 < m_A(m_A-2)(5m_A-2)^{-1}$.

In what follows we prove (140). The argument is similar to that used to prove tightness in Section 3.2, but since the estimates must be uniform in time, a more careful analysis is required. By (62) and (139),

$$\|\hat{X}^n\|^*_{T_n}, \|\hat{Y}^n\|^*_{T_n}, \|\hat{Z}^n\|^*_{T_n}, \|\Psi^n\|^*_{T_n} \le c_{11}(1+T)^m(1+\|\hat{W}^n\|^*_T) =: \xi(n,T).$$

Clearly

$$P(\vartheta_n \le T) \le \sum_{i\sim j} P\Big(\sup_{t\le T} \Lambda^n_{ij}(t) \ge b_1\Big).$$

Let $i\sim j$ be fixed. Recall that $b_1/2 \ge b_0$. As a result of (56)–(60),

$$P\Big(\sup_{t\le T} \Lambda^n_{ij}(t) \ge b_1\Big)$$
$$\le P(\exists\, s,r \in [\Theta_n, T_n], s \le r : \Lambda^n_{ij}(s) \le b_1/2, \tag{141}$$
$$\Lambda^n_{ij}(t) > 0, t \in [s,r], \Lambda^n_{ij}(r) \ge b_1).$$

Note that $b_1$ was chosen so that $1 + \max_{i\sim j}\sup_n \Lambda^n_{ij}(0) < b_1/2$. Note also that with $\Delta^n$ as in (115), (119) still holds. Arguing as in (118), using the fact that for $s \ge \Theta_n$, $h(t,\cdot) = h^0$, a constant function, with an appropriate constant $c_{12}$ one has

$$|\check{\Psi}^n_{ij}(t) - \check{\Psi}^n_{ij}(s)| \le c_{12}\xi(n,T)\|\hat{X}^n(t) - \hat{X}^n(s)\| + 1. \tag{142}$$

Also, $\bar{\Psi}^n - \psi^* = n^{-1/2}\hat{\Psi}^n$ and therefore $\|\Psi^n - \psi^*\|^*_{T_n} \le n^{-1/2}\xi(n,T)$. Thus for $s,t$ as in (141), recalling the notation (107) for $w_T(x,\delta)$, using (119),

$$n^{-1/2}\Delta^n_{ij}(s,t) \ge -2w_{c_{13}T_n}(\hat{S}^n_{ij}, c_{13}(t-s)) \tag{143}$$
$$- c_{13}\xi(n,T)(t-s) + \bar{c}n^{1/2}(t-s),$$

with $\bar{c} = \min \mu_{ij}\psi^*_{ij} > 0$. Combining (116), (111), (142), (143) and specializing to $t = r$,

$$
\begin{aligned}
b_1/2 &\le \Lambda^n_{ij}(t) - \Lambda^n_{ij}(s) \\
&\le 2w_{c_{13}T_n}(\hat{S}^n_{ij}, c_{13}(t-s)) \\
&\quad + c_{13}\xi(n,T)(t-s) - \bar{c}n^{1/2}(t-s) \\
&\quad + c_{12}\xi(n,T)\|\hat{X}^n(t) - \hat{X}^n(s)\|.
\end{aligned} \tag{144}
$$

Hence

$$P(\vartheta_n \le T) \le P(\Omega^n_1) + P(\Omega^n_2),$$

where

$$
\begin{aligned}
\Omega^n_1 &= \{\exists\, 0 \le s \le t \le T_n : t - s \le n^{-1/4},\ (144) \text{ holds}\}, \\
\Omega^n_2 &= \{\exists\, 0 \le s \le t \le T_n : t - s > n^{-1/4},\ (144) \text{ holds}\}.
\end{aligned}
$$

On $\Omega^n_1$,

$$b_1/2 \le 2w_{c_{13}T_n}(\hat{S}^n_{ij}, n^{-1/4}) + c_{13}\xi(n,T)n^{-1/4} + c_{12}\xi(n,T)\|\hat{X}^n_t - \hat{X}^n_s\|,$$

where $t - s \le n^{-1/4}$. By Lemma 9 below, recalling that $b_1 > 12$,

$$P(2w_{c_{13}T}(\hat{S}^n_{ij}, n^{-1/4}) \ge b_1/6) \le c_{15}Tn^{1/4 - m_A/8}.$$

By (65), and the bound $b_1$ on $\hat{M}^n$, we have $\|\hat{X}^n(t) - \hat{X}^n(s)\| \le c\xi(n,T)(t-s) + c\|\hat{W}^n(t) - \hat{W}^n(s)\| + b_1(t-s)$. Arguing again by Lemma 9, choosing $\beta_0$ large enough,

$$
\begin{aligned}
&P(c_{12}\xi(n,T)\max\{\|\hat{X}^n_t - \hat{X}^n_s\| : 0 \le t - s \le n^{-1/4}, t \le T_n\} \ge b_1/6) \\
&\quad \le cT^{\beta_1}n^{1/4 - m_A/8}.
\end{aligned}
$$

Also, by (98) and Lemma 8,

$$P(\xi(n,T)n^{-1/4} \ge b_1/6) \le P(\|\hat{W}^n\|^*_T \ge cn^{1/4}) \le c_{16}(1+T)^{m_A/2}n^{-m_A/4}.$$

Hence

$$P(\Omega^n_1) \le c_{15}T^{\beta_1}n^{1/4 - m_A/8} + c_{16}(1+T)^{m_A/2}n^{-m_A/4}. \tag{145}$$

On $\Omega^n_2$,

$$
\begin{aligned}
n^{1/4} &\le \|\hat{S}^n_{ij}\|^*_{c_{13}T_n} + c_{13}\xi(n,T)T_n + c_{12}\xi(n,T)\|\hat{X}^n_t - \hat{X}^n_s\| \\
&\le \|\hat{S}^n_{ij}\|^*_{c_{13}T_n} + c_{17}(1+T)^{c_{17}}(1 + \|\hat{W}^n\|^*_T)^2.
\end{aligned}
$$

Hence by (134),

$$P(\Omega^n_2) \le c_{18}(1+T)^{c_{18}}n^{-m_A/8}. \tag{146}$$

Combining (145) and (146) we obtain (140). Finally, the last statement in the proposition follows since all estimates used in this proof do not depend on the function $h_n$: we have used the fact that for $t \geq \Theta_n$, $h(t,\cdot) = h^0$.

Part (ii). Using (65), (66) and Gronwall's inequality [just like in the derivation of (104)], one has for $q = p'(\{h_n\})$ and an appropriate constant $C_1$ not depending on $t$, $n$ and $\{h_n\}$

$$\|\hat{X}^n\|_t^* \leq C_1 e^{C_1 t}(1 + \|\hat{W}^n\|_t^* + |\hat{M}^n|_t^*).$$

Replacing the estimate (130) with the above and reviewing the proof above of part (i), one recovers the estimate (77) as claimed. □

Lemma 9. *Given $\beta_0 \geq 0$ there are constants $c_1, \beta_1$, independent of $n$ and $T$ such that*

$$P(w_T(\hat{A}_i^n, n^{-1/4}) \geq T^{-\beta_0}) \leq c_1 T^{\beta_1} n^{1/4 - m_A/8}, \qquad n \in \mathbb{N}, T \geq 1.$$

*A similar estimate holds for $\hat{S}_{ij}^n$ in place of $\hat{A}_i^n$.*

Proof. Fix $i$ and suppress it from the notation. By (23),

$$\hat{A}^n(t) - \hat{A}^n(s) = n^{-1/2}(A^n(t) - A^n(s)) - n^{-1/2}\lambda^n(t-s).$$

Thus

$$P(w_T(\hat{A}^n, n^{-1/4}) \geq T^{-\beta_0}) \leq P(A^n(T) > 2\lambda^n T) + P(\Omega_{0,+}^n) + P(\Omega_{0,-}^n), \tag{147}$$

where

$$\begin{aligned}\Omega_{0,+}^n = \{&A^n(T) \leq 2\lambda^n T, \exists\, s,t \in [0, 2\lambda^n T], 0 \leq t-s \leq n^{-1/4},\\ &A^n(t) - A^n(s) \geq n^{1/2}T^{-\beta_0} + \lambda^n(t-s)\},\\ \Omega_{0,-}^n = \{&A^n(T) \leq 2\lambda^n T, \exists\, s,t \in [0, 2\lambda^n T], 0 \leq t-s \leq n^{-1/4},\\ &A^n(t) - A^n(s) \leq -n^{1/2}T^{-\beta_0} + \lambda^n(t-s)\}.\end{aligned}$$

Recall that $E\check{U}(k) = 1$. Letting $\bar{U}^n(k) = U^n(k) - (\lambda^n)^{-1}$, by (6) we have $E\bar{U}^n(k) = 0$. Let $M_j^n = \sum_{k=1}^j \bar{U}^n(k)$ and note that it is a martingale. For a real-valued function $X$ on $\mathbb{Z}_+$ let $\operatorname{osc}(X,i,j) = \max\{|X(k) - X(l)| : k,l \in [i+1,j]\}$. By (7), using $\lambda^n \leq c_2 n$, denoting $\rho = c_2 n^{3/4}$, we have

$$\begin{aligned}P(\Omega_{0,+}^n) &\leq P\left(\exists\, j \leq 2\lambda^n T, \exists\, r \leq n^{-1/4}\lambda^n, \sum_{k=j+1}^{j+n^{1/2}T^{-\beta_0}+r} U^n(k) \leq (\lambda^n)^{-1} r\right)\\ &\leq P\left(\exists\, j \leq 2c_2 nT, \exists\, r \leq c_2 n^{3/4}, \sum_{k=j+1}^{j+n^{1/2}+r} \bar{U}^n(k) \leq -(\lambda^n)^{-1} n^{1/2} T^{-\beta_0}\right)\end{aligned}$$

$$\leq P(\exists\, j \leq 2c_2nT, \operatorname{osc}(M^n, j, j+2c_2n^{3/4}) \geq c_2^{-1}n^{-1/2}T^{-\beta_0})$$
$$\leq P(\exists\, j \in [0, 2c_2nT] \cap \{0, \rho, 2\rho, \ldots\},$$
$$\operatorname{osc}(M^n, j, j+\rho) \geq c_2^{-1}n^{-1/2}T^{-\beta_0}/3)$$
$$\leq 1 - (1 - P(|M^n|_\rho^* \geq c_2^{-1}n^{-1/2}T^{-\beta_0}/6))^{2n^{1/4}T}.$$

Burkholder's inequality shows that

$$E(|M^n|_\rho^*)^{m_A} \leq c_3 E(\rho|\bar{U}(1)|^2)^{m_A/2} \leq c_4 n^{3m_A/8}(\lambda^n)^{-m_A} \leq c_5 n^{-5m_A/8}.$$

Hence

$$P(\Omega_{0,+}^n) \leq 1 - (1 - c_6 n^{-m_A/8} T^{m_A\beta_0})^{2n^{1/4}T} \leq c_7 n^{1/4 - m_A/8} T^{1+m_A\beta_0},$$

for an appropriate constant $c_7$. By a similar argument one shows that a similar bound holds for $P(\Omega_{0,-}^n)$. As follows from Lemma 8,

$$P(A^n(T) > 2\lambda^n T) \leq cn^{-m_A/2}(1+T)^{-m_A/2}.$$

Hence by (147),

$$P(w_T(\hat{A}^n, n^{-1/4}) \geq 1) \leq c_1 n^{1/4 - m_A/8} T^{1+m_A\beta_0}$$

for an appropriate constant $c_1$ independent of $n$ and $T \geq 1$. Since the Poisson processes $S_{ij}^n$ are in particular renewal processes (with finite $m_A$th moment for interarrival times), a similar result holds for $\hat{S}_{ij}^n$. $\square$

## APPENDIX

*Derivation of* (30). By (24), (21), (23)

$$\hat{X}_i^n(t) = n^{-1/2}(X_i^n(t) - nx_i^*)$$
$$= \hat{X}_i^{0,n} + \hat{A}_i^n(t) - \sum_j \hat{S}_{ij}^n\left(\int_0^t \bar{\Psi}_{ij}^n(s)\,ds\right) - \hat{R}_i^n\left(\int_0^t \bar{Y}_i^n(s)\,ds\right)$$
$$+ n^{-1/2}\left[\lambda_i^n t - \sum_j n\mu_{ij}^n \int_0^t \bar{\Psi}_{ij}^n(s)\,ds - n\theta_i^n \int_0^t \bar{Y}_i^n(s)\,ds\right].$$

Substituting (31), (19), (14) and (26) in the above, we obtain

$$\hat{X}_i^n(t) = \hat{X}_i^{0,n} + r_i \hat{W}_i^n(t) - \sum_j \mu_{ij}^n \int_0^t \hat{\Psi}_{ij}^n(s)\,ds$$
$$- \theta_i^n \int_0^t \hat{Y}_i^n(s)\,ds + \hat{\lambda}_i^n t + n^{1/2}\left[\lambda_i - \sum_j \mu_{ij}^n \psi_{ij}^*\right] t.$$

Finally, according to (14), (17) and (19), the last term above can be written as $-\sum_j \hat{\mu}_{ij}^n \psi_{ij}^* t$. Equation (30) follows. $\square$

**Acknowledgments.** I thank the referees for valuable suggestions that led to much improvement in the exposition.

## REFERENCES

[1] Ata, B. and Kumar, S. (2005). Heavy traffic analysis of open processing networks with complete resource pooling: Asymptotic optimality of discrete review policies. *Ann. Appl. Probab.* **15** 331–391. MR2115046

[2] Atar, R. (2005). A diffusion model of scheduling control in queueing systems with many servers. *Ann. Appl. Probab.* **15** 820–852. MR2114991

[3] Atar, R., Mandelbaum, A. and Reiman, M. (2004). Scheduling a multi-class queue with many exponential servers: Asymptotic optimality in heavy-traffic. *Ann. Appl. Probab.* **14** 1084–1134. MR2071417

[4] Bell, S. L. and Williams, R. J. (2001). Dynamic scheduling of a system with two parallel servers in heavy traffic with resource pooling: Asymptotic optimality of a threshold policy. *Ann. Appl. Probab.* **11** 608–649. MR1865018

[5] Billingsley, P. (1999). *Convergence of Probability Measures*, 2nd ed. Wiley, New York. MR1700749

[6] Ethier, S. N. and Kurtz, T. G. (1986). *Markov Processes. Characterization and Convergence*. Wiley, New York. MR838085

[7] Gans, N., Koole, G. and Mandelbaum, A. (2003). Telephone call centers: Tutorial, review and research prospects. *Manufacturing and Service Operations Management* **5** 79–141.

[8] Halfin, S. and Whitt, W. (1981). Heavy-traffic limits for queues with many exponential servers. *Oper. Res.* **29** 567–588. MR629195

[9] Harrison, J. M. (2000). Brownian models of open processing networks: Canonical representation of workload. *Ann. Appl. Probab.* **10** 75–103. MR1765204

[10] Harrison, J. M. and López, M. J. (1999). Heavy traffic resource pooling in parallel-server systems. *Queueing Systems Theory Appl.* **33** 339–368. MR1742575

[11] Karatzas, I. and Shreve, S. E. (1991). *Brownian Motion and Stochastic Calculus*, 2nd ed. Springer, New York. MR1121940

[12] Mandelbaum, A. and Reiman, M. I. Private communication.

[13] Mandelbaum, A. and Stolyar, A. L. (2004). Scheduling flexible servers with convex delay costs: Heavy traffic optimality of the generalized $c\mu$ rule. *Oper. Res.* **52** 836–855. MR2104141

[14] Williams, R. J. (2000). On dynamic scheduling of a parallel server system with complete resource pooling. In *Analysis of Communication Networks*: *Call Centres*, *Traffic and Performance*. *Fields Inst. Commun.* **28** 49–71. Amer. Math. Soc., Providence, RI. MR1788708

Department of Electrical Engineering
Technion–Israel Institute of Technology
Haifa 32000
Israel
e-mail: atar@ee.technion.ac.il